\documentclass{amsart}

\usepackage{amsmath,amsthm,amsfonts,amscd,amssymb}
\usepackage[all]{xy}

\newtheorem{thm}{Theorem}[section]
\newtheorem{cor}[thm]{Corollary}
\newtheorem{lem}[thm]{Lemma}

\theoremstyle{definition}

\theoremstyle{remark}
\newtheorem{rem}{Remark}[section]

\begin{document}

\title{On similarity classes of well-rounded sublattices of $\mathbb Z^2$}
\author{Lenny Fukshansky}

\address{Department of Mathematics, Claremont McKenna College, 850 Columbia Avenue, Claremont, CA 91711-6420}
\email{lenny@cmc.edu}

\begin{abstract}
A lattice is called well-rounded if its minimal vectors span the corresponding Euclidean space. In this paper we study the similarity classes of well-rounded sublattices of ${\mathbb Z}^2$. We relate the set of all such similarity classes to a subset of primitive Pythagorean triples, and prove that it has the structure of a noncommutative infinitely generated monoid. We discuss the structure of a given similarity class, and define a zeta function corresponding to each similarity class. We relate it to Dedekind zeta of ${\mathbb Z}[i]$, and investigate the growth of some related Dirichlet series, which reflect on the distribution of well-rounded lattices. We also construct a sequence of similarity classes of well-rounded sublattices of ${\mathbb Z}^2$, which gives good circle packing density and converges to the hexagonal lattice as fast as possible with respect to a natural metric we define. Finally, we discuss distribution of similarity classes of well-rounded sublattices of ${\mathbb Z}^2$ in the set of similarity classes of all well-rounded lattices in ${\mathbb R}^2$.
\end{abstract}

\maketitle
\tableofcontents

\def\A{{\mathcal A}}
\def\B{{\mathcal B}}
\def\C{{\mathcal C}}
\def\D{{\mathcal D}}
\def\E{{\mathcal E}}
\def\F{{\mathcal F}}
\def\x{{\mathcal H}}
\def\I{{\mathcal I}}
\def\J{{\mathcal J}}
\def\K{{\mathcal K}}
\def\L{{\mathcal L}}
\def\Ll{{\mathfrak L}}
\def\M{{\mathcal M}}
\def\Mm{{\mathfrak M}}
\def\Pp{{\mathfrak P}}
\def\Aa{{\mathfrak A}}
\def\Ss{{\mathfrak S}}
\def\N{{\mathcal N}}
\def\PP{{\mathcal P}}
\def\R{{\mathcal R}}
\def\s{{\mathcal S}}
\def\V{{\mathcal V}}
\def\W{{\mathcal W}}
\def\X{{\mathcal X}}
\def\Y{{\mathcal Y}}
\def\H{{\mathcal H}}
\def\OO{{\mathcal O}}
\def\aaa{{\mathbb A}}
\def\cee{{\mathbb C}}
\def\Nn{{\mathbb N}}
\def\pee{{\mathbb P}}
\def\que{{\mathbb Q}}
\def\real{{\mathbb R}}
\def\zed{{\mathbb Z}}
\def\gmn{{\mathbb G_m^N}}
\def\qbar{{\overline{\mathbb Q}}}
\def\DL{{\underline{\Delta}}}
\def\DU{{\overline{\Delta}}}
\def\eps{{\varepsilon}}
\def\vek{{\varepsilon_k}}
\def\ahat{{\hat \alpha}}
\def\bhat{{\hat \beta}}
\def\gt{{\tilde \gamma}}
\def\h{{\tfrac12}}
\def\ba{{\boldsymbol a}}
\def\be{{\boldsymbol e}}
\def\bei{{\boldsymbol e_i}}
\def\bc{{\boldsymbol c}}
\def\bm{{\boldsymbol m}}
\def\bk{{\boldsymbol k}}
\def\bi{{\boldsymbol i}}
\def\bl{{\boldsymbol l}}
\def\bq{{\boldsymbol q}}
\def\bu{{\boldsymbol u}}
\def\bt{{\boldsymbol t}}
\def\bs{{\boldsymbol s}}
\def\bv{{\boldsymbol v}}
\def\bw{{\boldsymbol w}}
\def\bx{{\boldsymbol x}}
\def\bX{{\boldsymbol X}}
\def\bz{{\boldsymbol z}}
\def\bwy{{\boldsymbol y}}
\def\bg{{\boldsymbol g}}
\def\bY{{\boldsymbol Y}}
\def\bL{{\boldsymbol L}}
\def\baa{{\boldsymbol\alpha}}
\def\bb{{\boldsymbol\beta}}
\def\bet{{\boldsymbol\eta}}
\def\bxi{{\boldsymbol\xi}}
\def\bo{{\boldsymbol 0}}
\def\bol{{\boldsymbol 1}_L}
\def\ep{\varepsilon}
\def\p{\boldsymbol\varphi}
\def\q{\boldsymbol\psi}
\def\WR{\operatorname{WR}}
\def\rank{\operatorname{rank}}
\def\aut{\operatorname{Aut}}
\def\lcm{\operatorname{lcm}}
\def\sgn{\operatorname{sgn}}
\def\spn{\operatorname{span}}
\def\md{\operatorname{mod}}
\def\Norm{\operatorname{Norm}}
\def\dim{\operatorname{dim}}
\def\det{\operatorname{det}}
\def\Vol{\operatorname{Vol}}
\def\rk{\operatorname{rk}}
\def\md{\operatorname{mod}}
\def\sqp{\operatorname{sqp}}
\def\Aut{\operatorname{Aut}}
\def\GL{\operatorname{GL}}
\def\Sim{\operatorname{Sim}}

\section{Introduction and statement of results}

Let $N \geq 2$ be an integer, and let $\Lambda \subseteq \real^N$ be a lattice of full rank. Define the {\it minimum} of $\Lambda$ to be
$$|\Lambda| = \min_{\bx \in \Lambda \setminus \{\bo\}} \|\bx\|,$$
where $\|\ \|$ stands for the usual Euclidean norm on $\real^N$. Let
$$S(\Lambda) = \{ \bx \in \Lambda : \|\bx\| = |\Lambda| \}$$
be the set of {\it minimal vectors} of $\Lambda$. We say that $\Lambda$ is a {\it well-rounded} lattice (abbreviated WR) if $S(\Lambda)$ spans $\real^N$. WR lattices come up in a wide variety of different contexts, including discrete optimization (e.g. sphere packing, covering, and kissing number problems), coding theory, and the linear Diophantine problem of Frobenius, just to name a few. In particular, the classical discrete optimization problems on lattices can usually be reduced to WR lattices in every dimension. Distribution of unimodular WR lattices in $\real^N$ has been studied by C. McMullen in \cite{mcmullen}. Also, the distribution of full-rank WR sublattices of $\zed^2$ has been recently studied in \cite{me:wr}. The goal of this paper is to continue this investigation from a somewhat different perspective. In particular, in \cite{me:wr} the zeta function $\zeta_{\WR(\zed^2)}(s)$ of WR sublattices of $\zed^2$ has been introduced, and we proved that it is analytic in the half-plane $\Re(s)>1$ with a pole of order {\it at least two} at $s=1$. In Theorem \ref{zeta_sim} we establish that in fact the order of the pole is {\it exactly} two, where the notion of the order of the pole we use here is defined by (\ref{pole_def}) below. To obtain this result we study the structure of the set of similarity classes of WR sublattices of $\zed^2$ (Theorems \ref{pyth_group}, \ref{main_par}, and \ref{sim_struct}) and use it to provide a simple analytic description for the Dirichlet series corresponding to each such similarity class (Theorem \ref{main:2}). We then decompose $\zeta_{\WR(\zed^2)}(s)$ over similarity classes to prove Theorem \ref{zeta_sim}. We also discuss sphere packing density of similarity classes of WR sublattices of $\zed^2$ (Theorem \ref{app_hex} and Corollary \ref{app_cor}), as well as their distribution among {\it all} WR similarity classes in~$\real^2$ (Theorem \ref{wr_dense}). 

We start out with a few words of motivation for the problems we study here. The similarity classes of the integer lattice $\zed^2$ and the hexagonal lattice $\Lambda_h$ (defined in (\ref{hex}) below) are very special in dimension two: these are the only two strongly eutactic similarity classes in $\real^2$, and $\left< \Lambda_h \right>$ is the only strongly perfect similarity class (we define the notions of strong eutaxy and perfection at the end of Section~5, see (\ref{design}) in particular; also see \cite{martinet}, especially chapter 16, for a detailed discussion of strongly eutactic and strongly perfect lattices and their properties). The distribution of sublattices of $\Lambda_h$ is studied in \cite{sloane}, while the distribution of {\it all} sublattices of $\zed^2$ is well understood (see for instance \cite{me:wr}, especially (52) and the beginning of Section~8, for a discussion of this). Studying the distribution of WR sublattices of these lattices is arguably even more important, since the WR property is vital in lattice theory. The goal of \cite{me:wr} and the current paper is to carry out this investigation for~$\zed^2$. We now introduce necessary notation and describe our results in more details.

Recall that two lattices $\Lambda_1, \Lambda_2 \subseteq \real^N$ of rank $N$ are said to be {\it similar} if there exists a matrix $A$ in $O_N(\real)$, the group of $N \times N$ real orthogonal matrices, and a real constant $\alpha$ such that $\Lambda_1 = \alpha A \Lambda_2$. This is an equivalence relation, which we will denote by writing $\Lambda_1 \sim \Lambda_2$, and the equivalence classes of lattices under this relation in $\real^N$ are called {\it similarity classes}. The distribution of sublattices of $\zed^N$ among similarity classes has been investigated by W. M. Schmidt in \cite{schmidt:sim}. 

The first trivial observation we can make is that WR property is preserved under similarity. In other words, if two full-rank lattices $\Lambda_1, \Lambda_2 \subseteq \real^N$ are similar, say $\Lambda_1 = \alpha A \Lambda_2$ for some $\alpha \in \real$ and $A \in O_N(\real)$, then
$$\det(\Lambda_1) = |\alpha|^N \det(\Lambda_2),\ \ |\Lambda_1| = |\alpha| |\Lambda_2|,$$
and $\Lambda_1$ is WR if and only if $\Lambda_2$ is WR. Therefore we can talk about similarity classes of well-rounded lattices in $\real^N$. From now on we will write $\WR(\Omega)$ for the set of all full-rank WR sublattices of a lattice $\Omega$; we will concentrate on $\WR(\zed^N)$, so let us write $\D_N$ and $\Mm_N$ for the sets of determinant and squared minima values, respectively, of lattices from $\WR(\zed^N)$. We will also write $\C_N$ for the set of all similarity classes of lattices in $\WR(\zed^N)$: this is a slight abuse of notation, since elements of $\C_N$ are really nonempty intersections of similarity classes of lattices in $\real^N$ with $\WR(\zed^N)$, as indicated in (\ref{sim_int}) below when $N=2$.
\bigskip

In this paper we study the case $N=2$. It is known that for every $\Lambda \in \WR(\zed^2)$ the set $S(\Lambda)$ has cardinality 4, and contains a {\it minimal basis} for $\Lambda$, which is unique up to $\pm$ signs and reordering (see Lemma~3.2 of \cite{me:wr}). For each $q \in \zed_{>0}$, define
\begin{equation}
\label{sq}
\s_q = \left\{ \frac{p}{q} \in \que \cap \left( \frac{\sqrt{3}}{2}, 1 \right)\ :\ \gcd(p,q) = 1, \sqrt{q^2 - p^2} \in \zed \right\}, 
\end{equation}
and let
\begin{equation}
\label{S}
\s = \left( \bigcup_{q \in \zed_{>0}} \s_q \right) \cup \{1\},
\end{equation}
where 1 is also thought of as $p/q$ with $p=q=1$. It is easy to see that the union in (\ref{S}) is disjoint, and each $\s_q$ is a subset of the set of Farey fractions of order $q$ in the interval $\left( \frac{\sqrt{3}}{2}, 1 \right)$. In Section 2 we show that the similarity classes of lattices in $\WR(\zed^2)$ are in bijective correspondence with fractions $p/q \in \s$. From now on, for each $p/q \in \s$, we will write $C(p,q)$ for the corresponding similarity class in $\C_2$, the set of all similarity classes of lattices in $\WR(\zed^2)$; a formal definition of $C(p,q)$ is given by (\ref{cpq_def}). The class $C(1,1)$ plays a special role: it is precisely the similarity class of all {\it orthogonal} well-rounded lattices, i.e. lattices of the form $\left( \begin{matrix} a&-b \\ b&a \end{matrix} \right) \zed^2$ for some $a,b \in \zed$. The set $\C_2$ has interesting algebraic and combinatorial structure. It is not difficult to notice that the set $\s$, which parametrizes $\C_2$, is in bijective correspondence with the set of primitive Pythagorean triples whose shortest leg is less than half of the hypothenuse. In Section 2 we explore this connection in details and use it to prove the following result.

\begin{thm} \label{pyth_group} The set $\C_2$ of similarity classes of lattices in $\WR(\zed^2)$ has the algebraic structure of an infinitely generated free non-commutative monoid with the class $C(1,1)$ of orthogonal well-rounded lattices serving as identity. As a combinatorial object, $\C_2$ has the structure of a regular rooted infinite tree, where each vertex has infinite degree, which is precisely the Cayley digraph of this monoid.
\end{thm}

\begin{rem} If $G$ is a monoid with a generating set $X$, then we define its {\it Cayley digraph} to be a directed graph with vertices corresponding to the elements of $G$, and with a directed edge between vertices $g$ and $h$ if $h=gx$ for some $x \in X$ (see for instance \cite{rotman} for details and related terminology).
\end{rem}

\noindent
We explicitly construct the monoid and the corresponding tree structure for $\C_2$ in Section 2. Notice that due to Theorem \ref{pyth_group} it makes sense to think of $\C_2$ as the moduli space of lattices in $\WR(\zed^2)$. 
\bigskip

In Section 3, we discuss a more explicit parametrization of $\C_2$, which allows to see the structure of each similarity class $C(p,q)$. It turns out that, although most well-rounded lattices are not orthogonal, all similarity classes in $\C_2$ can be parametrized by a subset of lattices from $C(1,1)$. More precisely, let us define a subset of $\zed^2$
\begin{eqnarray}
\label{A_par}
\A & = & \Big\{ (a,b) \in \zed^2 : 0<b<a,\ \gcd(a,b)=1,\ 2\nmid (a+b),\nonumber \\ 
& &\ \ \ \text{and either } b < a < \sqrt{3} b,\ \text{or } (2+\sqrt{3}) b < a  \Big\},
\end{eqnarray}
and consider the corresponding subset of $C(1,1)$
\begin{equation}
\label{orth_set}
C'(1,1) = \left\{ \left( \begin{matrix} a&-b \\ b&a \end{matrix} \right) \zed^2 \in C(1,1) : (a,b) \in \A \right\}.
\end{equation}
In Section 3 we prove the following theorem.

\begin{thm} \label{main_par} For each $p/q \in \s$, there exists a unique lattice $$\Omega = \left( \begin{matrix} a&-b \\ b&a \end{matrix} \right) \zed^2 \in C'(1,1),$$
where $a,b \in \zed_{>0}$ are given by
\begin{equation}
\label{pq_cond}
p=\max \{a^2-b^2, 2ab \}, \text{ and } q=a^2+b^2,
\end{equation}
such that $\Lambda \in C(p,q)$ if and only if
\begin{equation}
\label{lbx0}
\Lambda = \spn_{\zed} \left\{ \bx,\ \left( \begin{matrix} \frac{\sqrt{q^2-p^2}}{q}&-\frac{p}{q} \\ \frac{p}{q}&\frac{\sqrt{q^2-p^2}}{q} \end{matrix} \right) \bx \right\}
\end{equation}
for some $\bx \in \Omega$. Moreover, every lattice in the set $C'(1,1)$ parametrizes some similarity class $C(p,q)$ with $p,q$ as in (\ref{pq_cond}) in this way.
\end{thm}

\noindent
An easy consequence of Theorem \ref{main_par} is the existence of a lattice in each similarity class $C(p,q)$ which, in a sense to be described below, generates $C(p,q)$. First let us recall that given a full-rank lattice $\Lambda$ in $\real^2$, its Epstein zeta function is defined by
$$E_{\Lambda}(s) = \sideset{}{'}\sum_{\bx \in \Lambda} \|\bx\|^{-2s},$$
where $s \in \cee$, and $'$ indicates that the sum is taken over all $\bx \in \left( \Lambda / \{\pm1\} \right) \setminus \{\bo\}$. For each such $\Lambda$, this Dirichlet series is known to converge for all $s$ with $\Re(s) > 1$. Moreover, $E_{\Lambda}(s)$ has analytic continuation to $\cee$ except for a simple pole at $s=1$. For more information on $E_{\Lambda}(s)$ and its properties see \cite{sarnak}. In Section~3 we also prove the following theorem.

\begin{thm} \label{sim_struct} Let $C(p,q) \in \C_2$. There exists a lattice $\Lambda_{p,q} \in C(p,q)$, satisfying the following properties.

\begin{trivlist}
\item (1)\ $|\Lambda_{p,q}| = \min \{ |\Lambda| : \Lambda \in C(p,q) \} = \sqrt{q}$,
\item (2)\ $\det(\Lambda_{p,q}) = \min \{ \det(\Lambda) : \Lambda \in C(p,q) \} = p$,
\item(3)\ The norm form of $\Lambda_{p,q}$ with respect to its minimal basis is
$$Q_{p,q}(x,y) = qx^2+2xy\sqrt{q^2-p^2}+qy^2,$$
\item (4)\ For each $\Lambda \in C(p,q)$ there exists $U \in O_2(\real)$ such that $\Lambda = \sqrt{\frac{\det(\Lambda)}{p}}\ U \Lambda_{p,q}$; the quadratic form $\left( \frac{\det(\Lambda)}{p} \right) Q_{p,q}(x,y)$ is therefore the norm form for $\Lambda$ with respect to its minimal basis,
\item(5)\ The Epstein zeta function of any lattice $\Lambda \in C(p,q)$ is of the form
$$E_{\Lambda}(s) = \left( \frac{p}{\det(\Lambda)} \right)^s \sideset{}{'}\sum_{(x,y) \in \zed^2} \frac{1}{Q_{p,q}(x,y)^s},$$
and so $\Lambda_{p,q}$ maximizes $E_{\Lambda}(s)$ on $C(p,q)$ for each real value of $s > 1$.
\end{trivlist}
We call $\Lambda_{p,q}$ a {\it minimal} lattice of its similarity class $C(p,q)$; it is unique up to a rational rotation.
\end{thm}

Lattices $\Lambda_{p,q}$ also determine zeta-functions of corresponding similarity classes $C(p,q)$. Namely, with each $C(p,q) \in \C_2$ we can now associate two Dirichlet series, which incorporate information about the determinants and the minima of lattices in this similarity class, respectively. Specifically, define
$$Z^d_{p,q}(s) = \sum_{\Lambda \in C(p,q)} (\det(\Lambda))^{-s},\ \ \ \ \ Z^m_{p,q}(s) = \sum_{\Lambda \in C(p,q)} |\Lambda|^{-2s},$$
where $s \in \cee$. Our next goal is to investigate the properties of $Z^d_{p,q}(s)$ and $Z^m_{p,q}(s)$ for each $C(p,q) \in \C_2$, which we do by relating them to the Epstein zeta function of the lattice $\Omega$ parametrizing $C(p,q)$, as in Theorem \ref{main_par}. 
\smallskip

We will also write $\zeta_K(s)$ for the Dedekind zeta function of a number field $K$. It is known to be analytic for all $s \in \cee$ with $\Re(s) > 1-1/d$, where $d=[K:\que]$, except for a simple pole at $s=1$. For more information on properties of $\zeta_K(s)$ see~\cite{lang}. There is a standard relation between Dedekind zeta of imaginary quadratic fields and lattices of rank two, a special case of which we exploit here; see \cite{zagier} for more details.
\smallskip

\begin{thm} \label{main:2}  For each $C(p,q) \in \C_2$,
\begin{equation}
\label{zeta_r}
Z^d_{p,q}(s) = \frac{1}{p^s} \zeta_{\que(i)}(s)= \frac{1}{(\det(\Lambda_{p,q}))^s} \zeta_{\que(i)}(s),
\end{equation}
and
\begin{equation}
\label{zeta_r1}
Z^m_{p,q}(s) = \frac{1}{q^s} \zeta_{\que(i)}(s) = \frac{1}{|\Lambda_{p,q}|^{2s}} \zeta_{\que(i)}(s).
\end{equation}
\end{thm}

\noindent
We prove Theorem \ref{main:2} in Section 3, as well. Notice in particular that $C(1,1)$, the similarity class of all lattices coming from ideals in $\zed[i]$, has $Z^d_{1,1}(s) = Z^m_{1,1}(s) = \zeta_{\que[i]}(s)$, since $\Lambda_{1,1}=\zed^2$. This fact is also discussed in \cite{me:wr}.
\smallskip

In \cite{me:wr} we studied basic properties of the zeta-function of {\it all} well-rounded lattices
$$\zeta_{\WR(\zed^2)}(s) = \sum_{\Lambda \in \WR(\zed^2)} (\det(\Lambda))^{-s}.$$
It also makes sense to define
$$\zeta^m_{\WR(\zed^2)}(s) = \sum_{\Lambda \in \WR(\zed^2)} |\Lambda|^{-2s}.$$
These two Dirichlet series carry information about the distribution of lattices in $\WR(\zed^2)$ with respect to their determinant and minima values. For each similarity class $C(p,q) \in \C_2$, let us call $p$ its {\it determinant weight} and $q$ its {\it minima weight}. Theorem \ref{main:2} immediately implies that
\begin{eqnarray}
\label{wrd_zeta}
\zeta_{\WR(\zed^2)}(s) & = & \sum_{C(p,q) \in \C_2} \sum_{\Lambda \in C(p,q)} (\det(\Lambda))^{-s} \nonumber \\
& = & \sum_{C(p,q) \in \C_2} Z^d_{p,q}(s) = \zeta_{\que(i)}(s) \mathop{\sum_{p\ :\ p/q \in \s}}_{\text{for some } q \in \zed_{>0}}  \frac{a_p}{p^s},
\end{eqnarray}
and similarly
\begin{eqnarray}
\label{wrm_zeta}
\zeta^m_{\WR(\zed^2)}(s) = \zeta_{\que(i)}(s) \mathop{\sum_{q\ :\ p/q \in \s}}_{\text{for some } p \in \zed_{>0}}  \frac{b_q}{q^s},
\end{eqnarray}
where $a_p$ is the number of similarity classes in $\C_2$ with determinant weight $p$, and $b_q$ is the number of similarity classes in $\C_2$ with minima weight $q$; notice that $b_q = |\s_q|$, where $\s_q$ is as in (\ref{sq}). In fact, let us write
\begin{equation}
\label{weight_en}
W_d(s) = \sum_{C(p,q) \in \C_2} \frac{1}{p^s} = \mathop{\sum_{p\ :\ p/q \in \s}}_{\text{for some } q \in \zed_{>0}}  \frac{a_p}{p^s},
\end{equation}
and
\begin{equation}
\label{weight_enm}
W_m(s) = \sum_{C(p,q) \in \C_2} \frac{1}{q^s}  = \mathop{\sum_{q\ :\ p/q \in \s}}_{\text{for some } p \in \zed_{>0}}  \frac{b_q}{q^s}.
\end{equation}
We will call $W_d(s)$ and $W_m(s)$ {\it determinant} and {\it minima weight enumerators}, respectively.  Therefore the question of distribution of lattices in $\WR(\zed^2)$ is linked to understanding the basic analytic properties of $W_d(s)$ and $W_m(s)$. In Section~4 we use an approach different from that of \cite{me:wr} to prove the following result.

\begin{thm} \label{zeta_sim} Let the notation be as above, then $W_d(s)$ and $W_m(s)$ both have simple poles at $s=1$ and are analytic for all $s \in \cee$ with $\Re(s) > 1$. Therefore $\zeta_{\WR(\zed^2)}(s)$ and $\zeta^m_{\WR(\zed^2)}(s)$ both have poles of order two at $s=1$ and are analytic for all $s \in \cee$ with $\Re(s) > 1$.
\end{thm}

\noindent
We should point out that we are using the notion of a pole here {\it not} in a sense that would imply the existence of an analytic continuation, but only to reflect on the growth of the coefficients. More precisely, for a Dirichlet series $\sum_{n=1}^{\infty} c_n n^{-s}$, we say that it has a {\it pole of order} $\mu$ at $s=s_0$, where $\mu$ and $s_0$ are positive real numbers, if
\begin{equation}
\label{pole_def}
0 < \lim_{s \rightarrow s_0^+} |s-s_0|^{\mu} \sum_{n=1}^{\infty} |c_n n^{-s}| < \infty.
\end{equation}

\noindent
Notice that Theorem \ref{zeta_sim} in particular improves slightly on the result of Theorem 1.5 of \cite{me:wr}. The approach we use in Section~4 to prove Theorem \ref{zeta_sim} uses bounds on coefficients of weight enumerators $W_d(s)$ and $W_m(s)$ by coefficients of Dirichlet series associated with the set of primitive Pythagorean triples, which have Euler product expansions.
\bigskip

A standard object of lattice theory is a sphere packing associated with a lattice, and a classical problem is to determine the optimal packing density among lattices in a given dimension (see \cite{conway}). This problem has been solved in dimension two; in fact, it is not difficult to show that maximization of packing density can be restricted to WR lattices. Here we will discuss the circle packing density corresponding to lattices in $\WR(\zed^2)$, investigating how "close" can one come to the optimal packing density in dimension two with such lattices. For these purposes, let us write $\left< \Lambda \right>$ for the similarity class of any lattice $\Lambda$ in $\real^2$, so that
\begin{equation}
\label{def_sim}
\left< \Lambda \right> = \{ \alpha U \Lambda : \alpha \in \real_{>0},\ U \in O_2(\real) \}.
\end{equation}
Then for each $p/q \in \s$, 
\begin{equation}
\label{sim_int}
C(p,q) = \left< \Lambda_{p,q} \right> \cap \WR(\zed^2).
\end{equation}
For a lattice $\Lambda$ in $\real^2$ define
\begin{equation}
\label{theta_def}
\theta(\Lambda) = \min \left\{ \arcsin \left( \frac{| \bx^t \bwy |}{\|\bx\| \|\bwy\|} \right) : \bx,\bwy \text{ is a shortest basis for } \Lambda \right\}.
\end{equation}
By a {\it shortest basis} $\bx,\bwy$ of $\Lambda$ we mean here that $\bx$ is a minimal vector of $\Lambda$, and $\bwy$ is a vector of smallest Euclidean norm such that $\bx,\bwy$ is a basis for $\Lambda$. By a well known lemma of Gauss, $\theta(\Lambda) \in \left[ \frac{\pi}{3}, \frac{\pi}{2} \right]$ (see \cite{me:wr}). It is easy to notice that $\theta(\Lambda)$ remains constant on $\left< \Lambda \right>$, so we can also write $\theta(\left< \Lambda \right>)$. If $\bx,\bwy$ is a shortest basis for $\Lambda$ with the angle between $\bx$ and $\bwy$ equal to $\theta(\Lambda)$, then 
$$\det(\Lambda) = \|\bx\| \|\bwy\| \sin \theta(\Lambda),$$
and so if $\Lambda$ is well-rounded, then $\|\bx\| = \|\bwy\| = |\Lambda|$, and so
\begin{equation}
\label{wr_theta}
\det(\Lambda) = |\Lambda|^2 \sin \theta(\Lambda).
\end{equation}
It is easy to see that two well-rounded lattices $\Lambda_1,\Lambda_2 \subseteq \real^2$ are similar if and only if $\theta(\Lambda_1)=\theta(\Lambda_2)$, i.e. if and only if
$$\sin \theta(\Lambda_1) = \sin \theta(\Lambda_2) \in \left[ \frac{\sqrt{3}}{2}, 1 \right],$$
and so similarity classes of well-rounded lattices in $\real^2$ are indexed by real numbers in the interval $\left[ \frac{\sqrt{3}}{2}, 1 \right]$. Let $\Sim(\real^2)$ be the set of all similarity classes of well-rounded lattices in $\real^2$, and for every two $\left< \Lambda_1 \right>, \left< \Lambda_2 \right> \in \Sim(\real^2)$ define
\begin{equation}
\label{d_s}
d_s(\Lambda_1,\Lambda_2) = \left| \sin \theta(\Lambda_1) - \sin \theta(\Lambda_2)\right|.
\end{equation}
It is easy to see that $d_s$ is a metric on $\Sim(\real^2)$. If $\Lambda$ is a well-rounded lattice in $\real^2$, then the density of circle packing given by $\Lambda$ is
\begin{equation}
\label{pack1}
\delta(\Lambda) = \frac{\pi |\Lambda|^2}{4 \det(\Lambda)} = \frac{\pi}{4 \sin \theta(\Lambda)},
\end{equation}
by (\ref{wr_theta}), and so it depends not on the particular lattice $\Lambda$, but on its similarity class $\left< \Lambda \right>$. Moreover, (\ref{pack1}) implies that the smaller is $\sin \theta(\Lambda)$ the bigger is $\delta(\Lambda)$. Indeed, it is a well known fact that the similarity class $\left< \Lambda_h \right>$ gives the optimal circle packing in dimension two, where
\begin{equation}
\label{hex}
\Lambda_h = \left( \begin{matrix} 1&\frac{1}{2} \\ 0&\frac{\sqrt{3}}{2} \end{matrix} \right) \zed^2
\end{equation}
is the two-dimensional hexagonal lattice, and $\sin \theta(\Lambda_h) = \frac{\sqrt{3}}{2}$. The lattice $\Lambda_h$ also has the largest minimum among all lattices in $\real^2$ with the same determinant, and minimizes Epstein zeta function for all real values of $s > 1$ (see \cite{cassels}). However, $\left< \Lambda_h \right> \cap \WR(\zed^2) = \emptyset$. How well, with respect to the metric $d_s$ on $\Sim(\real^2)$, can we approximate the similarity class $\left< \Lambda_h \right>$ with similarity classes of the form $\left< \Lambda_{p,q} \right>$, i.e. with similarity classes that have a nonempty intersection with the set $\WR(\zed^2)$? This question is especially interesting since, in contrast to the two-dimensional situation, the three-dimensional counterpart of $\Lambda_h$, the face-centered cubic (fcc) lattice which maximizes sphere packing density in~$\real^3$, is in $\WR(\zed^3)$. Our next result addresses this question.

\begin{thm} \label{app_hex} There exists an infinite sequence of similarity classes $\left< \Lambda_{p_k,q_k} \right>$ such that
$$\left< \Lambda_{p_k,q_k} \right> \longrightarrow \left< \Lambda_h \right>, \text{ as } k \rightarrow \infty,$$
with respect to the metric $d_s$ on $\Sim(\real^2)$. The rate of this convergence can be expressed by  
\begin{equation}
\label{aphx}
\frac{1}{3\sqrt{3}\ q_k} < d_s(\Lambda_h,\Lambda_{p_k,q_k}) < \frac{1}{2\sqrt{3}\ q_k},
\end{equation}
where $q_k = O(14^k)$ as $k \rightarrow \infty$. Moreover, the inequality (\ref{aphx}) is sharp in the sense that
\begin{equation}
\label{aphx2}
\frac{1}{3\sqrt{3}\ q} < d_s(\Lambda_h,\Lambda_{p,q}),
\end{equation}
for every similarity class of the form $\left< \Lambda_{p,q} \right> \neq \left< \Lambda_{1,1} \right>$. For the similarity class of orthogonal well-rounded lattices $\left< \Lambda_{1,1} \right> = \left< \zed^2 \right>$, we clearly have $d_s(\Lambda_h,\zed^2) = \frac{2-\sqrt{3}}{2}$.
\end{thm}

\begin{cor} \label{app_cor} Each similarity class $\left< \Lambda_{p_k,q_k} \right>$ of Theorem \ref{app_hex} gives circle packing density $\delta_{p_k,q_k}$ such that
\begin{equation}
\label{pack2.1}
\delta(\Lambda_h) \left( \frac{1}{1 + \frac{1}{723 \times (13.928)^{k-1}}} \right) < \delta_{p_k,q_k} < \delta(\Lambda_h) \left( \frac{1}{1 + \frac{0.92}{723 \times (13.947)^{k-1}}} \right),
\end{equation}
where $\delta(\Lambda_h) = \frac{\pi}{\sqrt{12}} = 0.9069...$ is the circle packing density of $\Lambda_h$.
\end{cor}

\noindent
We prove Theorem \ref{app_hex} and Corollary \ref{app_cor} in Section~5. Notice that a well-rounded lattice in $\real^2$ has a rational basis, i.e. a basis consisting of vectors with rational coordinates, if and only if it belongs to a similarity class $\left< \Lambda_{p,q} \right>$ for some $p,q$. Therefore results of Theorem \ref{app_hex} and Corollary \ref{app_cor} can be interpreted as statements on best approximation to $\Lambda_h$ (and hence best circle packing) by well-rounded lattices in $\real^2$ with rational bases. As we will see in Section~5, this just comes down to finding best approximations to $\frac{\sqrt{3}}{2}$ by fractions $\frac{p}{q}$ where $(p,\sqrt{q^2-p^2},q)$ is a primitive Pythagorean triple with $\sqrt{q^2-p^2} \leq q/2$. In fact, a similar approximation result holds for all WR lattices in $\real^2$, not just $\Lambda_h$.

\begin{thm} \label{wr_dense} The similarity classes of WR sublattices of $\zed^2$ are dense in the set of all similarity classes of WR lattices in $\real^2$, in other words the set $\{ \left< \Lambda_{p,q} \right> : p/q \in \s\}$ is dense in $\Sim(\real^2)$ with respect to the metric $d_s$. Moreover, for every $\Lambda \in \Sim(\real^2)$, there exist infinitely many non-similar lattices $\Lambda_{p,q} \in \WR(\zed^2)$ such that 
\begin{equation}
\label{wr_dense_bound}
d_s( \Lambda, \Lambda_{p,q} ) \leq \frac{2 \sqrt{2}}{q}.
\end{equation}
\end{thm}

\noindent
We derive Theorem \ref{wr_dense} in Section 6 as an easy corollary of a theorem of Hlawka on Diophantine approximation with quotients of Pythagorean triples, and discuss equidistribution of $\{ \left< \Lambda_{p,q} \right> : p/q \in \s\}$ in $\Sim(\real^2)$. As a side remark in Section 6, we also use Hlawka's result to approximate points on a rational ellipse by rational points on the same ellipse. Notice that Theorem \ref{wr_dense} does not include Theorem \ref{app_hex} as a special case, since the approximating constants in Theorem \ref{app_hex} are sharper and the proof is constructive unlike that of Theorem \ref{wr_dense}. We are now ready to proceed.
\bigskip

\section{Parametrization by Pythagorean triples}

Notice that if a lattice $\Lambda \in \WR(\zed^2)$, then $\cos \theta(\Lambda), \sin \theta(\Lambda) \in \que_{>0}$, where $\theta(\Lambda)$ is defined in (\ref{theta_def}), and  therefore we can index similarity classes of lattices in $\WR(\zed^2)$ by fractions $p/q \in \s$, where $\s$ is as in (\ref{S}), so for each such $p/q$ the corresponding similarity class $C(p,q) \in \C_2$ is a set of the form
\begin{equation}
\label{cpq_def}
C(p,q) = \left\{ \Lambda \in \WR(\zed^2)\ :\ \sin \theta(\Lambda) = \frac{p}{q} \right\}.
\end{equation}
For each $p/q \in \s$, define $t = \sqrt{q^2-p^2} \in \zed$. Then it is easy to notice that 
$$0 \leq t < \frac{q}{2} < \frac{\sqrt{3} q}{2} < p \leq q,$$
and $t^2+p^2=q^2$ with $\gcd(t,p,q)=1$. In other words, the set $\s$, and therefore the set $\C_2$ of similarity classes of lattices in $\WR(\zed^2)$, is in bijective correspondence with the set of primitive Pythagorean triples with the shortest leg being less than half of the hypothenuse. 
\smallskip

Let
$$\Pp = \{ (x,y,z)\ :\ x,y,z \in \zed_{>0},\ 2|y,\ \gcd(x,y,z)=1,\ x^2+y^2=z^2 \}$$ 
be the set of all primitive Pythagorean triples, and let
$$\PP = \{ (x,y,z) \in \Pp\ :\ \min \{x,y\} < z/2 \} \cup \{(1,0,1)\}.$$
Notice that we include $(1,0,1)$ in $\PP$, although it is traditionally not included in $\Pp$. Then $p/q \in \s$ if and only if either $(t,p,q) \in \PP$ or $(p,t,q) \in \PP$. In other words, elements of $\PP$ can be used to enumerate similarity classes of lattices in $\WR(\zed^2)$. We will use this approach to provide a convenient combinatorial description of elements of $\C_2$. Define matrices
\begin{equation}
\label{ABC}
A = \left( \begin{matrix} 1&-2&2 \\ 2&-1&2 \\ 2&-2&3 \end{matrix} \right),\ B = \left( \begin{matrix} 1&2&2 \\ 2&1&2 \\ 2&2&3 \end{matrix} \right),\ C = \left( \begin{matrix} -1&2&2 \\ -2&1&2 \\ -2&2&3 \end{matrix} \right) \in \GL_3(\zed),
\end{equation}
and let $G = \left< I_3, A, B, C \right>$ be the non-commutative monoid generated by $A,B,C$ with the $3 \times 3$ identity matrix $I_3$. Let us think of elements of $\Pp$ as vectors in $\zed^3$, and for each $M \in G$ define the corresponding linear transformations
\begin{equation}
\label{action}
M(x,y,z) = M\left( \begin{matrix} x\\y\\z \end{matrix} \right).
\end{equation}
It is a well known fact that for every $(x,y,z) \in \Pp$, $A(x,y,z), B(x,y,z), C(x,y,z) \in \Pp$. Moreover, every $(x,y,z) \in \Pp$ can be obtained in a unique way by applying a sequence of linear transformations $A,B,C$ to $(3,4,5)$, the smallest triple in $\Pp$ (this construction is attributed to Barning \cite{barning}; also see \cite{alperin}, \cite{romik}). This means that (\ref{action}) defines a free action of $G$ on the set $\Pp$ of primitive Pythagorean triples by left multiplication. The set $\Pp$ has the structure of an infinite rooted ternary tree with respect to this action, as described in \cite{alperin}; this in particular implies that $G$ is a free monoid. In fact, this tree (see Figure 1 below) is precisely the Cayley digraph of $G$ with respect to the generating set $\{A,B,C\}$.

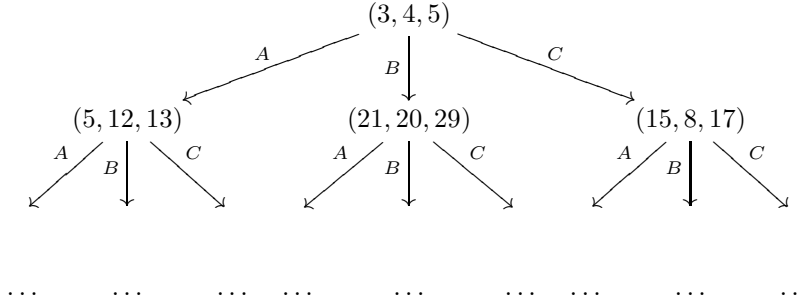
\begin{figure}[!ht]
\begin{center}
 \SelectTips{cm}{}
\[
  \xymatrix@C=7pt{
  &&&& (3,4,5) \ar_{A}[llld] \ar_{B}[d] \ar^{C}[rrrd]&
  \\
  & (5,12,13) \ar_{A}[ld] \ar_{B}[d] \ar^{C}[rd] &&
  & (21,20,29) \ar_{A}[ld] \ar_{B}[d] \ar^{C}[rd] &&
  & (15,8,17) \ar_{A}[ld] \ar_{B}[d] \ar^{C}[rd] &&
  \\
  & & & & & & & & 
  \\
  \dots& \dots& \dots& \dots& \dots& \dots& \dots& \dots& \dots& }
\]
\caption{Ternary tree representation for $\Pp$}
\end{center}
\end{figure}

We can  extend this construction by considering the set $\Pp' = \Pp \cup \{ (1, 0, 1) \}$ (compare with \cite{cass}). It is easy to notice that $A(1,0,1)=B(1,0,1)=(3,4,5)$ and $C(1,0,1)=(1,0,1)$, and hence every $(x,y,z) \in \Pp$ can be obtained by applying a sequence of linear transformations $A,B,C$ to $(1,0,1)$. Such a sequence is no longer unique, hence action of $G$ does not extend to $\Pp'$, however a {\it shortest} such sequence is unique up to multiplication on the left by either $BA^{-1}$ or $AB^{-1}$.

Notice that $\PP \subset \Pp'$. Let 
\begin{eqnarray}
\label{H_group}
H & = & \Big\{ A^2NB, A^kB, ABNB, C^2NB, C^kB, CBNB \nonumber \\
& &\ \ (AC)^kA^2NB, (AC)^kABNB, (AC)^kAB, \nonumber \\
& &\ \ (CA)^kC^2NB, (CA)^kCBNB, (CA)^kCB\ :\ N \in G, k \in \zed_{>0}  \Big\}.
\end{eqnarray}
It is clear that $H$ is a subsemigroup and $H' = H \cup \{ I_3 \}$ is a submonoid of $G$. Let us define the image of $\Pp'$ under $G$ to be
$$G\Pp' = \{ M(x,y,z)\ :\ M \in G,\ (x,y,z) \in \Pp' \},$$
and similarly for the images $G\PP$, $H\Pp'$, and $H\PP$.

\begin{lem} \label{pyth} $H\Pp' = H\PP = \PP \setminus \{(1,0,1)\}$.
\end{lem}

\proof
First we will prove that $H\Pp' \subseteq \PP$. It is clear that $(1,0,1) \notin H\Pp'$. Let $M \in H$, then there exists some $N \in G$ such that one of the following is true:
\begin{trivlist}
\item (1)\ $M = A^2 N$,
\item (2)\ $M = AB N$,
\item (3)\ $M = C^2 N$,
\item (4)\ $M = CB N$,
\item (5)\ $M = (AC)^kA^2 N$, where $k \in \zed_{>0}$,
\item (6)\ $M = (AC)^kAB N$, where $k \in \zed_{>0}$,
\item (7)\ $M = (CA)^kC^2 N$, where $k \in \zed_{>0}$,
\item (8)\ $M = (CA)^kCB N$, where $k \in \zed_{>0}$.
\end{trivlist}
Let $(x,y,z) \in \Pp'$, and write $(x',y',z') = N(x,y,z)$, where $N$ is as above. Then, in case (1) 
$$M(x,y,z) = A^2(x',y',z') = \left( \begin{matrix} x'-4y'+4z' \\ 4x'-7y'+8z' \\ 4x'-8y'+9z' \end{matrix} \right),$$
where
$$\frac{1}{2} (4x'-8y'+9z') = 2x'-4y'+\frac{9}{2}z' > x'-4y'+4z',$$
hence $M(x,y,z) \in \PP$. In case (2)
$$M(x,y,z) = AB(x',y',z') = \left( \begin{matrix} x'+4y'+4z' \\ 4x'+7y'+8z' \\ 4x'+8y'+9z' \end{matrix} \right),$$
where
$$\frac{1}{2} (4x'+8y'+9z') = 2x'+4y'+\frac{9}{2}z' > x'+4y'+4z',$$
hence $M(x,y,z) \in \PP$. In case (3)
$$M(x,y,z) = C^2(x',y',z') = \left( \begin{matrix} -7x'+4y'+8z' \\ -4x'+y'+4z' \\ -8x'+4y'+9z' \end{matrix} \right),$$
where
$$\frac{1}{2} (-8x'+4y'+9z') = -4x'+2y'+\frac{9}{2}z' > -4x'+y'+4z',$$
hence $M(x,y,z) \in \PP$. In case (4)
$$M(x,y,z) = CB(x',y',z') = \left( \begin{matrix} 7x'+4y'+8z' \\ 4x'+y'+4z' \\ 8x'+4y'+9z' \end{matrix} \right),$$
where
$$\frac{1}{2} (8x'+4y'+9z') = 4x'+2y'+\frac{9}{2}z' > 4x'+y'+4z',$$
hence $M(x,y,z) \in \PP$. 

For cases (5) and (6), let 
$$(x_2,y_2,z_2) = M(x,y,z) = (AC)^k(x_1,y_1,z_1),$$ 
where
$$(x_1,y_1,z_1) = A^2(x',y',z') = \left( \begin{matrix} x'-4y'+4z' \\ 4x'-7y'+8z' \\ 4x'-8y'+9z' \end{matrix} \right)$$
in case (5), and
$$(x_1,y_1,z_1) = AB(x',y',z') = \left( \begin{matrix} x'+4y'+4z' \\ 4x'+7y'+8z' \\ 4x'+8y'+9z' \end{matrix} \right)$$
in case (6). It is not difficult to notice that $x_2 = \min\{x_2,y_2\}$, and
$$\frac{z_2}{2} = x_2 + \left( \frac{z_1}{2} - x_1 \right).$$
Therefore $(x_2,y_2,z_2) \in \PP$ if and only if $x_1 \leq z_1/2$, which is true in both cases, (5) and (6). On the other hand,
$$B(x',y',z') = \left( \begin{matrix} x'+2y'+2z' \\ 2x'+y'+2z' \\ 2x'+2y'+3z' \end{matrix} \right),$$
and 
$$C(x',y',z') = \left( \begin{matrix} -x'+2y'+2z' \\ -2x'+y'+2z' \\ -2x'+2y'+3z' \end{matrix} \right),$$
which implies that $(AC)^kB(x',y',z'),\ (AC)^kC(x',y',z') \notin \PP$ for any $(x',y',z')$.

For cases (7) and (8), let 
$$(x_2,y_2,z_2) = M(x,y,z) = (CA)^k(x_1,y_1,z_1),$$ 
where
$$(x_1,y_1,z_1) = C^2(x',y',z') = \left( \begin{matrix} -7x'+4y'+8z' \\ -4x'+y'+4z' \\ -8x'+4y'+9z' \end{matrix} \right)$$
in case (7), and
$$(x_1,y_1,z_1) = CB(x',y',z') = \left( \begin{matrix} 7x'+4y'+8z' \\ 4x'+y'+4z' \\ 8x'+4y'+9z' \end{matrix} \right)$$
in case (8). It is not difficult to notice that $y_2 = \min\{x_2,y_2\}$, and
$$\frac{z_2}{2} = y_2 + \left( \frac{z_1}{2} - y_1 \right).$$
Therefore $(x_2,y_2,z_2) \in \PP$ if and only if $y_1 \leq z_1/2$, which is true in both cases, (7) and (8). On the other hand,
$$B(x',y',z') = \left( \begin{matrix} x'+2y'+2z' \\ 2x'+y'+2z' \\ 2x'+2y'+3z' \end{matrix} \right),$$
and 
$$A(x',y',z') = \left( \begin{matrix} x'-2y'+2z' \\ 2x'-y'+2z' \\ 2x'-2y'+3z' \end{matrix} \right),$$
which implies that $(CA)^kB(x',y',z'),\ (CA)^kA(x',y',z') \notin \PP$ for any $(x',y',z')$. We have shown that $H\PP \subseteq H\Pp' \subseteq \PP \setminus \{(1,0,1)\}$.

To finish the proof of the lemma, we will show that $\PP \setminus \{(1,0,1)\} \subseteq H\PP$. Notice that it is in fact sufficient to show that for each $(x,y,z) \in \PP \setminus \{(1,0,1)\}$ there exists $M \in H$ such that $(x,y,z) = M(1,0,1)$. We know that there exists $N \in G$ such that $(x,y,z) = N(3,4,5)$, and so $(x,y,z) = NB(1,0,1)$. First notice that $N$ cannot be of the form $BN'$ for some $N' \in G$. Indeed, suppose it is, then
$$(x,y,z) = B(x',y',z') = \left( \begin{matrix} x'+2y'+2z' \\ 2x'+y'+2z' \\ 2x'+2y'+3z' \end{matrix} \right),$$
where $(x',y',z') = N'(x,y,z) \in \Pp'$, but
$$\frac{1}{2} (2x'+2y'+3z') = x'+y'+\frac{3}{2}z' < \min \{ x'+2y'+2z', 2x'+y'+2z' \},$$
which contradicts the fact that $(x,y,z) \in \PP$. Similarly, from the arguments above it follows that $N$ cannot be of the form $(AC)^kBN'$, $(AC)^kCN'$, $(CA)^kBN'$, or $(CA)^kAN'$. The only options left are those described in cases (1) - (8) above, which means that $M=NB \in H$. Therefore $\PP \setminus \{(1,0,1)\} \subseteq H\PP$, which completes the proof.
\endproof

\begin{thm} \label{bij} $H'$ is a free infinitely generated monoid, which acts freely on the set $\PP$ by left multiplication. With respect to this action, $\PP$ has the structure of a regular rooted infinite tree, where each vertex has infinite degree (see Figure 2 below); this is precisely the Cayley digraph of $H'$.
\end{thm}

\proof
$H'$ is a submonoid of $G$, which is a free monoid, hence $H'$ must also be free by the Nielsen-Schreier theorem (see for instance \cite{rotman}). To see that $H'$ is infinitely generated, consider for instance the set $\{ AB^k : k \in \zed_{>0} \}$ of elements of $H'$. Since $A,B^k \notin H'$ for any $k \in \zed_{>0}$, it is clear that no finite subset of $H'$ can generate all of the elements of the form  $AB^k$: if this was possible, there would have to be relations between elements of $H'$, contradicting the fact that it is free. Therefore $H'$ must be infinitely generated, and so its Cayley digraph is  a regular rooted infinite tree, where each vertex has infinite degree, and the root corresponds to $I_3$.

By Lemma \ref{pyth} we know that $H'\PP = \PP$. Moreover, we know that for each $(x,y,z) \in \PP$ there exists a unique element $N \in G$ such that $N(3,4,5) = (x,y,z)$, and hence $NB$ is the unique element in $H$ such that $NB(1,0,1) = (x,y,z)$. This means that $H'$ acts freely on $\PP$. Then we can identify $(1,0,1) \in \PP$ with $I_3 \in H'$, and each $(x,y,z) \in \PP$ with the corresponding unique $NB \in H'$ such that $NB(1,0,1) = (x,y,z)$, which means that with respect to the action of $H'$ the set $\PP$ has the structure of the Cayley digraph of $H'$ with respect to an appropriate generating set.
\endproof

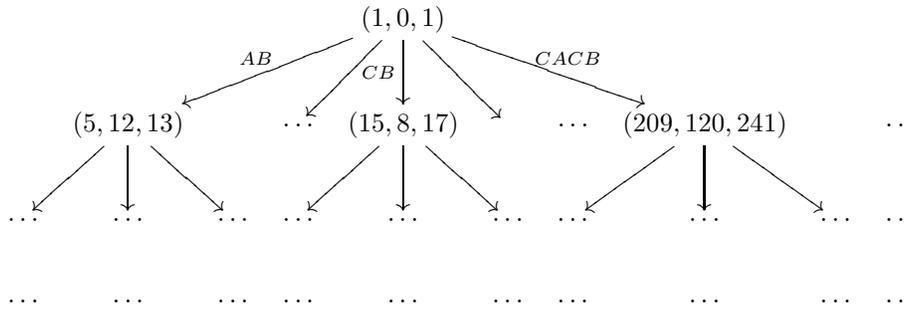
\begin{figure}[!ht]
\begin{center}
 \SelectTips{cm}{}
\[
  \xymatrix@C=7pt{
  &&&& (1,0,1) \ar_{AB}[llld] \ar[ld] \ar_{CB}[d] \ar[rd] \ar^{CACB}[rrrd] &
  \\
 & (5,12,13) \ar[ld] \ar[d] \ar[rd]  && \dots
  & (15,8,17) \ar[ld] \ar[d] \ar[rd] && \dots
  & (209,120,241) \ar[ld] \ar[d] \ar[rd] && \dots
  \\
  \dots& \dots& \dots& \dots& \dots& \dots& \dots& \dots& \dots& \dots
  \\
  \dots& \dots& \dots& \dots& \dots& \dots& \dots& \dots& \dots& \dots}
\]
\caption{Infinite-degree tree representation for $\PP$}
\end{center}
\end{figure}

\begin{cor} \label{group} The set $\C_2$ of similarity classes of lattices in $\WR(\zed^2)$ has the structure of a non-commutative free infinitely generated monoid. Specifically, it is isomorphic to $H'$.
\end{cor}

\proof
We will identify $\C_2$ with $H'$ in the following way. From Theorem \ref{bij} we know that there exists a bijection $\varphi : \PP \rightarrow H'$, given by
$$\varphi(x,y,z) = M, \text{ such that } M(1,0,1) = (x,y,z),$$
for each $(x,y,z) \in \PP$ with $\varphi^{-1} : H' \rightarrow \PP$ defined by
$$\varphi^{-1}(M) = M(1,0,1),$$
for each $M \in H'$.

On the other hand, there also exists a bijection $\psi: \C_2 \rightarrow \PP$, given by
\[ \psi(C(p,q)) = \left\{ \begin{array}{ll}
(\sqrt{q^2-p^2}, p, q) & \mbox{if $2 | p$} \\
(p, \sqrt{q^2-p^2}, q) & \mbox{if $2 \nmid p$,}
\end{array}
\right. \]
for each $C(p,q) \in \C_2$ with $\psi^{-1} : \PP \rightarrow \C_2$ defined by
$$\psi^{-1}(x,y,z) = C(p,q), \text{ where } p=\max\{x,y\},\ q=z,$$
for each $(x,y,z) \in \PP$. Therefore we have bijections $\varphi \psi : \C_2 \rightarrow H'$ and $(\varphi \psi)^{-1} = \psi^{-1} \varphi^{-1} : H' \rightarrow \C_2$.

We can now define a binary operation $*$ on $\C_2$ as follows: for every $C(p_1,q_1)$ and $C(p_2,q_2)$ in $\C_2$, let
\begin{equation}
\label{oper_*}
C(p_1,q_1) * C(p_2,q_2) = \psi^{-1} \left( \varphi \psi( C(p_1,q_1)) \varphi \psi( C(p_2,q_2)) (1,0,1) \right).
\end{equation}
It is easy to see that $\C_2$ is a free non-commutative monoid with respect to $*$, which is isomorphic to $H'$ via the monoid isomorphism $\varphi \psi : \C_2 \rightarrow H'$, and $(\varphi \psi)^{-1}(I_3) = C(1,1) \in \C_2$ is the identity. Hence the tree in Figure 2 is the Cayley digraph of $\C_2$ with respect to an appropriate generating set. This completes the proof.
\endproof

Now Theorem \ref{pyth_group} follows by combining Theorem \ref{bij} with Corollary \ref{group}.
\bigskip

\section{Similarity classes and corresponding zeta functions}

In this section we discuss the structure of similarity classes $C(p,q)$, as well as the properties of associated zeta functions. Our first goal is to prove Theorem \ref{main_par}. For each $p/q \in \s$, define
\begin{eqnarray*}
\Mm_2(p,q) = \Big\{ \bx \in \left( \zed^2 / \{\pm 1\} \right) \setminus \{\bo\} & : & x_1\sqrt{q^2 - p^2} \equiv x_2 p\ (\md q), \\
& & x_2 \sqrt{q^2 - p^2} \equiv -x_1 p\ (\md q) \Big\}.
\end{eqnarray*}
Notice that $\bx \in \Mm_2(p,q)$ if and only if
\begin{equation}
\label{lbx}
\Lambda(\bx) := \spn_{\zed} \left\{ \bx,\ \left( \begin{matrix} \frac{\sqrt{q^2-p^2}}{q}&-\frac{p}{q} \\ \frac{p}{q}&\frac{\sqrt{q^2-p^2}}{q} \end{matrix} \right) \bx \right\} \in C(p,q).
\end{equation}
Hence lattices in the similarity class $C(p,q)$ are in bijective correspondence with points in $\Mm_2(p,q)$.

\begin{lem} \label{one_less} Let $p/q \in \s$. The congruence relations
\begin{equation}
\label{cong1}
x_1 \sqrt{q^2 - p^2} \equiv x_2 p\ (\md q)
\end{equation}
and
\begin{equation}
\label{cong2}
x_2 \sqrt{q^2 - p^2} \equiv -x_1 p\ (\md q)
\end{equation}
are equivalent, meaning that
$$\Mm_2(p,q) = \left\{ \bx \in \left( \zed^2 / \{\pm 1\} \right) \setminus \{\bo\}\ :\ x_1 \sqrt{q^2 - p^2} \equiv x_2 p\ (\md q) \right\}.$$
\end{lem}

\proof
Recall that $\gcd(p,q)=1$. Also $q^2-p^2 = (p-q)(p+q)$, and
$$\gcd(q-p,q) = \gcd(q+p,q) = \gcd(p,q) = 1,$$
therefore $\gcd(q^2-p^2,q) = \gcd(\sqrt{q^2-p^2},q) = 1$. Hence
$$x_1 \sqrt{q^2 - p^2} \equiv x_2 p\ (\md q)$$
if and only if
$$x_2 p \sqrt{q^2 - p^2} \equiv x_1 (q^2 - p^2) \equiv -x_1 p^2\ (\md q),$$
which happens if and only if $-x_1 p \equiv x_2 \sqrt{q^2 - p^2}\ (\md q)$.
\endproof

For each $p/q \in \s$, define $c(p,q)$ to be the unique integer such that $0 \leq c(p,q) \leq q-1$ and
\begin{equation}
\label{cpq}
c(p,q) p \equiv \sqrt{q^2-p^2}\ (\md q).
\end{equation}
Then, by Lemma \ref{one_less}, for every $\bx \in \Mm_2(p,q)$ we have $x_2 \equiv c(p,q) x_1\ (\md q)$, meaning that $x_2 = c(p,q) x_1 + yq$ for some $y \in \zed$. In other words, $\Mm_2(p,q)$ can be presented as
\begin{equation}
\label{m2_rep}
\Mm_2(p,q) = \left\{ \left( \begin{matrix} x\\c(p,q)x+qy \end{matrix} \right)\ :\ \left( \begin{matrix} x\\y \end{matrix} \right) \in \left( \zed^2 / \{\pm 1\} \right) \setminus \{\bo\} \right\}.
\end{equation}
Define
\begin{equation}
\label{opq}
\Omega(p,q) = \left( \begin{matrix} 1&0 \\ c(p,q)&q \end{matrix} \right) \zed^2,
\end{equation}
so that $\Mm_2(p,q) = \left( \Omega(p,q) / \{\pm 1\} \right) \setminus \{\bo\}$.

\begin{lem} \label{opq_min} $\Omega(p,q) \in C(1,1)$, and $|\Omega(p,q)|^2 = \det(\Omega(p,q)) = q$. In fact, each $\Omega(p,q)$ is in the set $C'(1,1)$ as defined by (\ref{orth_set}). Moreover, every lattice $\left( \begin{matrix} a&-b \\ b&a \end{matrix} \right) \zed^2$ in the set $C'(1,1)$ is of the form $\Omega(p,q)$ for $p,q$ satisfying (\ref{pq_cond}).
\end{lem}

\proof
First fix a lattice $\Omega(p,q)$. It is a well known fact that there exist unique relatively prime $a > b \in \zed_{>0}$ of different parity such that either $p=a^2-b^2$ or $p=2ab$, and $q=a^2+b^2$ (this is the standard parametrization of primitive Pythagorean triples, see for instance \cite{silverman}). The fact that $\frac{\sqrt{3}}{2} q < p \leq q$ ensures that $(a,b) \in \A$. Then the lattice
$$\Omega = \left( \begin{matrix} a&-b \\ b&a \end{matrix} \right) \zed^2$$
is in $C'(1,1)$, and $|\Omega|^2 = \det(\Omega) = q$. We will now show that $\Omega(p,q) = \Omega$. Since $\gcd(a,b)=1$, there exist $g_1,g_2 \in \zed$ such that
\begin{equation}
\label{p0}
g_1a-g_2b=1.
\end{equation}
Let $\gamma = g_1b+g_2a$, and notice that
$$\left( \begin{matrix} a&-b \\ b&a \end{matrix} \right) \left( \begin{matrix} g_1&b \\ g_2&a \end{matrix} \right) = \left( \begin{matrix} 1&0 \\ \gamma&q \end{matrix} \right),$$
and $\det \left( \begin{matrix} a&-b \\ b&a \end{matrix} \right) = \det \left( \begin{matrix} 1&0 \\ \gamma&q \end{matrix} \right) = q$, so
\begin{equation}
\label{p1}
\Omega = \left( \begin{matrix} 1&0 \\ \gamma&q \end{matrix} \right) \zed^2.
\end{equation}
Notice that $\gamma^2+1$ is divisible by $q$. Indeed, (\ref{p0}) implies that $g_1 = \frac{g_2b+1}{a}$, and so 
\begin{eqnarray}
\label{p2}
\gamma^2+1 & = & (g_1b+g_2a)^2 + 1 = g_2^2a^2 + 2g_2b(g_2b+1) + \frac{b^2}{a^2} (g_2b+1)^2 + 1 \nonumber \\
& = & q \left( \frac{g_2q+2g_2b+1}{a^2} \right) = q (g_1^2+g_2^2).
\end{eqnarray}
Moreover, we can ensure that $0 \leq \gamma \leq q-1$ by replacing $\gamma$ with $\gamma + qm$ for some $m \in \zed$, if necessary: it is easy to see that $\gamma^2+1$ will still be divisible by $q$, and (\ref{p1}) will still hold. We will now show that $\gamma = c(p,q)$. Notice that
\begin{equation}
\label{gamma_q}
F_{\gamma,q}(x,y) = \left( \frac{\gamma^2+1}{q} \right) x^2 + 2\gamma xy + qy^2
\end{equation}
is an integral binary quadratic form with discriminant -4, hence it is equivalent to
$$G(x,y) = x^2+y^2,$$
since the class number of -4 is one. In fact, it is easy to verify that
$$F_{\gamma,q}(x,y) = G(g_1x-by,\ g_2x-ay),\ G(x,y) = F_{\gamma,q}(ax-by,\ g_2x-g_1y).$$
Let us write $t = a^2-b^2$, so either $p=t$ or $p = \sqrt{q^2-t^2}$, then
$$q = G(a,b) = F_{\gamma,q}(t,k) = \left( \frac{\gamma^2+1}{q} \right) t^2 + 2\gamma tk + qk^2,$$
where $k=g_2a-g_1b$. Therefore
$$\gamma^2t^2 + 2\gamma tk + q^2k^2 = q^2 - t^2,$$
meaning that
\begin{equation}
\label{p3}
\gamma^2t^2 \equiv q^2-t^2\ (\md q).
\end{equation}
Notice that $\gcd(\gamma,q) = 1$, since $q | (\gamma^2+1)$, so $\gcd(\gamma,q)$ must divide 1. Therefore, if $p=t$, then (\ref{p3}) implies that
\begin{equation}
\label{p4}
\gamma p \equiv \sqrt{q^2-p^2}\ (\md q).
\end{equation}
If, on the other hand, $p = \sqrt{q^2-t^2}$, then (\ref{p3}) implies that $\gamma \sqrt{q^2-p^2} \equiv p\ (\md q)$, meaning that $\gamma^2 \sqrt{q^2-p^2} \equiv \gamma p\ (\md q)$, but on the other hand $\gamma^2 \equiv -1\ (\md q)$, and so
$$-\sqrt{q^2-p^2} \equiv \gamma p\ (\md q).$$
By Lemma \ref{one_less}, this last congruence is equivalent to (\ref{p4}). We conclude that $0 \leq \gamma \leq q-1$, and $\gamma$ satisfies (\ref{p4}), which means that $\gamma = c(p,q)$, and so $\Omega = \Omega(p,q)$.
\smallskip

In the opposite direction, assume that
$$\Omega = \left( \begin{matrix} a&-b \\ b&a \end{matrix} \right) \zed^2 \in C'(1,1),$$
and define $q = a^2+b^2$, $p=\max \{a^2-b^2, 2ab \}$. Then the fact that $(a,b) \in \A$ ensures that $\frac{\sqrt{3}}{2} q < p \leq q$, i.e. $p/q \in \s$. It is not difficult to notice that for $p/q \in \s$, $p=a^2-b^2$ if and only if $a > (2+\sqrt{3}) b$, and $p=2ab$ if and only if $b < a < \sqrt{3} b$. The argument identical to the one above now shows that $\Omega = \Omega(p,q)$. This completes the proof.
\endproof

\begin{rem} It is not difficult to conclude from an argument very similar to the one in the proof of Lemma \ref{opq_min} that all binary integral quadratic forms of discriminant -4 are of the form $F_{\gamma,q}(x,y)$ as in (\ref{gamma_q}) for some odd positive integer $q$ which is not divisible by any prime of the form $4k+3$ and an integer $\gamma$ (positive or negative) such that $\gamma^2+1$ is divisible by $q$. This statement is essentially equivalent to the fact that there is a bijection between ideals of the form $(da+dbi)$ in $\zed[i]$ with $a>b>0$, $d>0$, and Pythagorean triples $(d^2(a^2-b^2),2d^2ab,d^2(a^2+b^2))$.
\end{rem}

\begin{proof}[Proof of Theorem \ref{main_par}]
Fix a similarity class $C(p,q) \in \C_2$ for some $p/q \in \s$. For each $\Lambda \in C(p,q)$, we have $\Lambda=\Lambda(\bx)$ as defined by (\ref{lbx}) for some $\bx = \left( \begin{matrix} x\\c(p,q)x+qy \end{matrix} \right)$ where $x,y \in \zed$, hence $\bx \in \Omega(p,q)$. On the other hand, for each $\bx \in \Omega(p,q)$, the corresponding lattice $\Lambda(\bx)$ is easily seen to be in $C(p,q)$. Combining this observation with Lemma \ref{opq_min} completes the proof of the theorem.
\end{proof}
\bigskip

\noindent
{\it Proof of Theorem \ref{sim_struct}.}
We first introduce the notion of a minimal lattice in each similarity class $C(p,q)$. Let $\bx(p,q) \in \Omega(p,q)$ be such that $\|\bx(p,q)\| = |\Omega(p,q)|$, and let $\Lambda(\bx(p,q))$ be defined by (\ref{lbx}); we will call this lattice a {\it minimal lattice} of the similarity class $C(p,q)$ and will denote it by $\Lambda_{p,q}$. By Lemma \ref{opq_min}, we have
\begin{equation}
\label{c_min}
|\Omega(p,q)| = \|\bx(p,q)\| = \sqrt{q}.
\end{equation}
On the other hand, since $\sqrt{3}/{2} < p/q \leq 1$, meaning that the angle between $\bx(p,q)$ and the other minimal basis vector given in (\ref{lbx}) is $\arcsin \left( p/q \right) \in \left( \pi/3, \pi/2 \right)$, a well known lemma of Gauss (see \cite{near:ort} or \cite{me:wr}) implies that $|\Lambda_{p,q}| = \|\bx(p,q)\|$, and so $|\Lambda_{p,q}| =\sqrt{q}$. Therefore
\begin{equation}
\label{lxc_det}
\det(\Lambda_{p,q}) = |\Lambda_{p,q}|^2\ \frac{p}{q} = p = \min\{ \det(\Lambda) : \Lambda \in C(p,q) \}.
\end{equation}
Moreover, a straight-forward computation shows that the norm form of $\Lambda_{p,q}$ with respect to its minimal basis is
\begin{equation}
\label{Q_pq}
Q_{p,q}(x,y) = (x,y) \left( \begin{matrix} q&\sqrt{q^2-p^2} \\ \sqrt{q^2-p^2}&q \end{matrix} \right) \left( \begin{matrix} x \\ y \end{matrix} \right) = qx^2 + 2xy \sqrt{q^2-p^2} + qy^2.
\end{equation}
Notice that the minimal lattice of a similarity class may not in general be unique, however it is unique up to a rational rotation, and so for our purposes it suffices to pick any one of them.
\smallskip

Next, let $\Lambda \in C(p,q)$, then $\Lambda \sim \Lambda_{p,q}$, and so there must exist $\alpha \in \real_{>0}$ and $U \in O_2(\real)$ such that $\Lambda = \alpha U \Lambda_{p,q}$. Then $\det(\Lambda) = \alpha^2 p$, and so $\alpha = \sqrt{\frac{\det(\Lambda)}{p}} > 1$. If we write $A$ and and $A_{p,q}$ for the minimal basis matrices of $\Lambda$ and $\Lambda_{p,q}$ respectively, then $A = \sqrt{\frac{\det(\Lambda)}{p}} U A_{p,q}$, and the norm form of $\Lambda$ with respect to this minimal basis is
$$Q_{\Lambda}(x,y) = (x,y) A^tA \left( \begin{matrix} x \\ y \end{matrix} \right) = \frac{\det(\Lambda)}{p} (x,y) A_{p,q}^tA_{p,q} \left( \begin{matrix} x \\ y \end{matrix} \right) = \frac{\det(\Lambda)}{p} Q_{p,q}(x,y).$$
Epstein zeta function of $\Lambda$ is therefore given by
$$E_{\Lambda}(s) = \sideset{}{'}\sum_{\bx \in \Lambda} \|\bx\|^{-2s} = \sideset{}{'}\sum_{(x,y) \in \zed^2} Q_{\Lambda}(x,y)^{-s} = \left( \frac{p}{\det(\Lambda)} \right)^s \sideset{}{'}\sum_{(x,y) \in \zed^2} Q_{p,q}(x,y)^{-s}.$$
Then (\ref{lxc_det}) implies that for every fixed real value of $s>1$, $E_{\Lambda}(s)$ achieves its maximum on $C(p,q)$ when $\Lambda=\Lambda_{p,q}$, and it does not achieve a minimum since there exist lattices in $C(p,q)$ with arbitrarily large determinants. This completes the proof of the theorem.
\boxed{ }
\bigskip

\begin{proof}[Proof of Theorem \ref{main:2}]
We now derive the properties of the Dirichlet series corresponding to each~$C(p,q)$. Fix a similarity class $C(p,q) \in \C_2$. By Theorem \ref{main_par}, each $\Lambda \in C(p,q)$ is of the form $\Lambda(\bx)$ for some $\bx \in \Omega(p,q)$. As in the proof of Theorem \ref{sim_struct} above, a well known lemma of Gauss (see \cite{near:ort} or \cite{me:wr}) implies that $|\Lambda(\bx)| = \|\bx\|$. Since also, by Lemma~3.2 of \cite{me:wr}, the set of minimal vectors of $\Lambda(\bx)$ is precisely 
$$\left\{ \pm \bx,\ \pm \left( \begin{matrix} \frac{\sqrt{q^2-p^2}}{q}&-\frac{p}{q} \\ \frac{p}{q}&\frac{\sqrt{q^2-p^2}}{q} \end{matrix} \right) \bx \right\},$$
it follows that $\Lambda(\bx_1) = \Lambda(\bx_2)$ if and only if $\bx_1 = \pm \bx_2$. Therefore
$$Z^m_{p,q}(s) = \sum_{\Lambda \in C(p,q)} |\Lambda|^{-2s} = \sideset{}{'}\sum_{\bx \in \Omega(p,q)} \|\bx\|^{-2s} = E_{\Omega(p,q)}(s).$$
Now Theorem \ref{main_par} readily implies that there exists $U \in O_2(\real)$ such that
$$\Omega(p,q) = U \left( \begin{matrix} \sqrt{q}&0 \\ 0&\sqrt{q} \end{matrix} \right) \zed^2,$$
which means that $E_{\Omega(p,q)}(s)$ is equal to the Epstein zeta function of $\left( \begin{matrix} \sqrt{q}&0 \\ 0&\sqrt{q} \end{matrix} \right) \zed^2$. Hence
$$Z^m_{p,q}(s) = E_{\Omega(p,q)}(s) = \frac{1}{q^s} \sideset{}{'}\sum_{\bx \in \zed^2} \|\bx\|^{-2s} = \frac{1}{q^s} \zeta_{\que(i)}(s),$$
which proves (\ref{zeta_r1}). Now recall that for each $\Lambda \in C(p,q)$,
$$\det(\Lambda) = |\Lambda|^2 \sin \theta(\Lambda) = |\Lambda|^2 \frac{p}{q}.$$
Then (\ref{zeta_r}) follows.
\end{proof}
\bigskip

Finally, we present a simple lemma, which is not related to the rest of this section, but is of some independent interest. Recall that the {\it dual} of a lattice $\Lambda$ in $\real^N$ is the lattice $\Lambda^*$, defined by
$$\Lambda^* = \{ \bx \in \real^N\ :\ \forall\ \bwy \in \Lambda,\ \bx^t \bwy \in \zed \}.$$

\begin{lem} \label{dual} Let $C(p,q) \in \C_2$, and suppose that $\Lambda \in C(p,q)$. Then $\det(\Lambda) \Lambda^* \in C(p,q)$.
\end{lem}

\proof
Let $\bx, \bwy$ be the minimal basis for $\Lambda$, and write $A = \left( \begin{matrix} x_1&y_1 \\ x_2&y_2 \end{matrix} \right)$ for the corresponding basis matrix. Then 
\begin{equation}
\label{d1}
(A^{-1})^t = \frac{1}{\det(\Lambda)} \left( \begin{matrix} y_2&-x_2 \\ -y_1&x_1 \end{matrix} \right)
\end{equation}
is a basis matrix for $\Lambda^*$ (see \cite{martinet}, p. 24). It is easy to notice that $\frac{1}{\det(\Lambda)} \left( \begin{matrix} y_2 \\ -y_1 \end{matrix} \right)$, $\frac{1}{\det(\Lambda)} \left( \begin{matrix} -x_2 \\ x_1 \end{matrix} \right)$ is therefore the minimal basis for $\Lambda^*$, and hence $\det(\Lambda) \Lambda^* \in \WR(\zed^2)$. Moreover, (\ref{d1}) implies that
$$\det(\Lambda) \Lambda^* = \left( \begin{matrix} 0&-1 \\ 1&0 \end{matrix} \right) \Lambda,$$
and so $\Lambda \sim \Lambda^*$, since $\left( \begin{matrix} 0&-1 \\ 1&0 \end{matrix} \right) \in O_2(\real)$.
\endproof
\bigskip

\section{Weight enumerators $W_d(s)$ and $W_m(s)$}

In this section we will discuss in more details some properties of the Dirichlet series $W_d(s)$ and $W_m(s)$ as defined in (\ref{weight_en}) and (\ref{weight_enm}), respectively, and will prove Theorem \ref{zeta_sim}. Recall that we write $a_p$ and $b_q$ for the coefficients of $W_d(s)$ and $W_m(s)$ respectively as defined in Section~1. The following formulas for $a_p$ and $b_q$ are immediate from Theorem \ref{main_par} and the definition of the set $\A$ in~(\ref{A_par}).

\begin{lem} \label{ap} For each $p$ such that $p/q \in \s$ for some $q \in \zed_{>0}$,
$$a_p = \left| \{ (m,n) \in \A : p= \max\{m^2-n^2,2mn\} \} \right|,$$
and for each $q$ such that $p/q \in \s$ for some $p \in \zed_{>0}$,
$$b_q = \left| \{ (m,n) \in \A : q=m^2+n^2 \} \right|.$$
\end{lem}

Notice that the expression for $a_p$ in Lemma \ref{ap} is similar in spirit to the function $\beta$ defined in \cite{me:wr}, in particular it can also be bounded in terms of Hooley's $\Delta$-function. On the other hand, we can obtain simple explicit bounds for $a_p$ and $b_q$ from our Pythagorean tree construction in Section~3. For each $p,q \in \zed_{>0}$, define $L(p)$ to be the number of primitive Pythagorean triples with a leg $p$, and $H(q)$ to be the number of primitive Pythagorean triples with the hypotenuse $q$. There are well known formulas for $L(p)$ and $H(q)$ (see \cite{beiler}, p. 116): if $p,q > 1$, then
\[ L(p) = \left\{ \begin{array}{ll}
0 & \mbox{if $p \equiv 2\ (\md 4)$} \\
2^{\omega(p)-1} & \mbox{otherwise,}
\end{array}
\right. \]
where $\omega(p)$ is the number of distinct prime divisors of $p$, and
\[ H(q) = \left\{ \begin{array}{ll}
0 & \mbox{if $2 | q$, or if $q$ has a prime factor $l \equiv 3\ (\md 4)$} \\
2^{\omega(q)-1} & \mbox{otherwise.}
\end{array}
\right. \]
For convenience, we also set $L(1) = H(1) = \frac{1}{2}$. It is clear that $a_p \leq L(p)$ and $b_q \leq H(q)$ when $p,q>1$, and $a_1=b_1=1$. One can ask how good are these bounds? We will now show that the correct order of magnitude of the bound for both, $a_p$ and $b_q$, in the sense that the corresponding Dirichlet series has the same behavior at $s=1$ as $W_d(s)$ and $W_m(s)$, is given by $H$ and not by $L$. Namely, define
$$\L(s) = \sum_{n=1}^{\infty} \frac{L(n)}{n^s},\ \ \ \H(s) = \sum_{n=1}^{\infty} \frac{H(n)}{n^s}.$$

\begin{lem} \label{hyp} $\H(s)$ has a simple pole at $s=1$ and is analytic for all $s \in \cee$ with $\Re(s) > 1$. Moreover when $\Re(s) > 1$, $\H(s)$ has an Euler product type expansion 
\begin{equation}
\label{euler1}
\H(s) = \frac{1}{2} \prod_{l \equiv 1 (\md 4)} \frac{l^s+1}{l^s-1},
\end{equation}
where the product is over primes $l$.
\end{lem}

\proof
Let us define
$$V_1 = \{ n \in \zed_{>0} : n \text{ is only divisible by primes which are } \equiv 1 (\md 4) \}.$$
Then notice that, as in the proof of Lemma 8.1 of \cite{me:wr},
\begin{eqnarray*}
2\H(s) = \sum_{n \in V_1} \frac{2^{\omega(n)}}{n^s} & = & \prod_{l \equiv 1 (\md 4)} \left( \sum_{k=0}^{\infty} 2^{\omega(l^k)} l^{-ks} \right) = \prod_{l \equiv 1 (\md 4)} \left( 1 + 2 \sum_{k=1}^{\infty} l^{-ks} \right) \\
& = & \prod_{l \equiv 1 (\md 4)} \left( \frac{2}{1-l^{-s}} - 1 \right) = \prod_{l \equiv 1 (\md 4)} \frac{l^s+1}{l^s-1},
\end{eqnarray*}
whenever this product is convergent, where $l$ is always prime. The fact that $\H(s)$ has a simple pole at $s=1$ and is analytic for all $s \in \cee$ with $\Re(s) > 1$ then follows immediately from Lemma 8.1 of \cite{me:wr}.
\endproof

\begin{proof}[Proof of Theorem \ref{zeta_sim}]
First of all notice that since $\frac{\sqrt{3}}{2} q \leq p \leq q$, we have
$$\left| \left( \frac{\sqrt{3}}{2} \right)^s \right| \sum_{C(p,q) \in \C_2} \left| \frac{1}{p^s} \right| \leq \sum_{C(p,q) \in \C_2} \left| \frac{1}{q^s} \right| \leq \sum_{C(p,q) \in \C_2} \left| \frac{1}{p^s} \right|,$$
meaning that $W_d(s)$ and $W_m(s)$ must have poles of the same order and the same half-plane of convergence. Since $b_q \leq H(q)$, Lemma \ref{hyp} implies that $W_m(s)$ has at most a simple pole at $s=1$, and is analytic when $\Re(s) > 1$. On the other hand, Theorem 1.5 of \cite{me:wr} implies that $\zeta_{\WR(\zed^2)}(s)$ has at least a pole of order two at $s=1$, meaning that, by (\ref{wrd_zeta}), $W_d(s)$ must have at most a simple pole at $s=1$. This means that both, $W_d(s)$ and $W_m(s)$, have simple poles at $s=1$ and are analytic when $\Re(s) > 1$, and therefore, by (\ref{wrd_zeta}) and (\ref{wrm_zeta}), $\zeta_{\WR(\zed^2)}(s)$ and $\zeta^m_{\WR(\zed^2)}(s)$ both have poles of order two at $s=1$ and are analytic when $\Re(s) > 1$. This completes the proof.
\end{proof}

\begin{rem}
\label{rem_sim}
Notice that Theorem \ref{zeta_sim} combined with Lemma \ref{hyp} implies that the Dirichlet series $\sum_{(x,y,z) \in \PP} \frac{1}{\max\{x,y\}^s}$ and $\sum_{(x,y,z) \in \Pp'} \frac{1}{\max\{x,y\}^s}$ have poles of the same order 1 at $s=1$. This fact could be roughly interpreted to mean that the sets $\PP$ and $\Pp'$ are comparable in size, i.e. that ``most'' primitive Pythagorean triples correspond to similarity classes of lattices from $\WR(\zed^2)$. In other words, the imposed condition that the shortest leg of a primitive Pythagorean triple is no longer than half of the hypothenuse is not particularly restrictive. Moreover, we can roughly think of $H(n)$ as a bound on the average orders of $a_n$ and $b_n$ for each $n \in \zed_{>0}$.
\end{rem}

On the other hand, we have the following.

\begin{lem} \label{legs} Let the notation be as above, then $\L(s)$ has a pole of order two at $s=1$ and is analytic for all $s \in \cee$ with $\Re(s) > 1$. Moreover when $\Re(s) > 1$, $\L(s)$ has an Euler product type expansion 
\begin{equation}
\label{euler}
\L(s) = \frac{1}{2} \left( \frac{4^s-2^s+2}{4^s-2^s} \right) \prod_{l \neq 2 \text{\ prime}} \frac{l^s+1}{l^s-1} = \frac{1}{2} \left( \frac{4^s-2^s+2}{4^s+2^s} \right) \frac{\zeta(s)^2}{\zeta(2s)},
\end{equation}
where $\zeta(s)$ is the Riemann zeta function.
\end{lem}

\proof
Let us consider the Dirichlet series $2 \L(s)$, then 
\begin{eqnarray}
\label{l1}
2 \L(s) & = & \sum_{n=1}^{\infty} \frac{2L(n)}{n^s} = \sum_{2 \nmid n} \frac{2^{\omega(n)}}{n^s} + \sum_{4|n} \frac{2^{\omega(n)}}{n^s} = \sum_{2 \nmid n} \frac{2^{\omega(n)}}{n^s} + \frac{1}{4^s} \sum_{n=1}^{\infty} \frac{2^{\omega(2n)}}{n^s} \nonumber \\
& = & \sum_{2 \nmid n} \frac{2^{\omega(n)}}{n^s} + \frac{1}{4^s} \left( 2 \sum_{2 \nmid n} \frac{2^{\omega(n)}}{n^s} + \sum_{2 | n} \frac{2^{\omega(n)}}{n^s} \right) \nonumber \\ 
& = & \left( 1 + \frac{2}{4^s} \right) \sum_{2 \nmid n} \frac{2^{\omega(n)}}{n^s} + \frac{1}{4^s}  \sum_{2 | n} \frac{2^{\omega(n)}}{n^s}.
\end{eqnarray}
On the other hand,
$$\sum_{2 | n} \frac{2^{\omega(n)}}{n^s} = \frac{1}{2^s} \sum_{n=1}^{\infty} \frac{2^{\omega(2n)}}{n^s} = \frac{1}{2^s} \left( 2 \sum_{2 \nmid n} \frac{2^{\omega(n)}}{n^s} + \sum_{2 | n} \frac{2^{\omega(n)}}{n^s} \right),$$
and so
\begin{equation}
\label{l2}
\sum_{2 | n} \frac{2^{\omega(n)}}{n^s} = \frac{2}{2^s-1} \sum_{2 \nmid n} \frac{2^{\omega(n)}}{n^s}.
\end{equation}
Combining (\ref{l1}) and (\ref{l2}), we obtain
\begin{equation}
\label{l3}
2 \L(s) = \left( \frac{4^s-2^s+2}{4^s-2^s} \right) \sum_{2 \nmid n} \frac{2^{\omega(n)}}{n^s}.
\end{equation}
Now define
\[ L_1(n) = \left\{ \begin{array}{ll}
0 & \mbox{if $2 | n$} \\
2^{\omega(n)} & \mbox{$2 \nmid n$.}
\end{array}
\right. \]
It is easy to see that $L_1(1) = 1$ and $L_1$ is multiplicative, i.e. if $\gcd(m,n)=1$ then $L_1(mn)=L_1(m)L_1(n)$. Therefore, by Theorem 286 of \cite{hardy},
\begin{eqnarray}
\label{l4}
\sum_{2 \nmid n} \frac{2^{\omega(n)}}{n^s} & = & \sum_{n=1}^{\infty} \frac{L_1(n)}{n^s} = \prod_{l \text{\ prime}} \left( \sum_{k=0}^{\infty} \frac{L_1(l^k)}{l^{ks}} \right) = \prod_{l \neq 2 \text{\ prime}} \left( 1 + 2 \sum_{k=1}^{\infty} l^{-ks} \right) \nonumber \\
& = & \prod_{l \neq 2 \text{\ prime}} \left( \frac{2}{1-l^{-s}} - 1 \right) = \prod_{l \neq 2 \text{\ prime}} \frac{l^s+1}{l^s-1},
\end{eqnarray}
when $\Re(s) > 1$. Moreover, by Theorem 301 of \cite{hardy},
\begin{equation}
\label{l5}
\frac{\zeta(s)^2}{\zeta(2s)} = \sum_{n=1}^{\infty} \frac{2^{\omega(n)}}{n^s} = \prod_{l \text{\ prime}} \frac{l^s+1}{l^s-1}.
\end{equation}
Now (\ref{euler}) follows by combining (\ref{l3}), (\ref{l4}), and (\ref{l5}). Moreover, $\zeta(s)^2/\zeta(2s)$ clearly has a pole of order two at $s=1$, and is analytic for all $s \in \cee$ with $\Re(s) > 1$. This completes the proof.
\endproof

\begin{rem}
\label{rem_sim1}
Since $\L(s)$ is the sum over {\it all} the legs of primitive Pythagorean triples, short and long, and it is easy to see that for each $(x,y,z) \in \Pp'$, $\max\{x,y\} \geq \frac{1}{\sqrt{2}} z$, Lemma \ref{legs} combined with Remark \ref{rem_sim} imply that $\sum_{(x,y,z) \in \Pp'} \frac{1}{\min\{x,y\}^s}$ must have a pole of order 2 at $s=1$.
\end{rem}

\bigskip

\section{Approximating the hexagonal lattice}

In this section we will talk about circle packing density corresponding to similarity classes of lattices in $\WR(\zed^2)$. Our goal is to prove Theorem \ref{app_hex}. We do it by first proving the following slightly more technical lemma, from which the theorem follows easily.

\begin{lem} \label{approx} Let $A,B,C$ be matrices as in (\ref{ABC}). For each $k \in \zed_{>0}$, let 
\begin{equation}
\label{ap_0}
(p_k,t_k,q_k) = (CA)^k CB(1,0,1) \in \PP,
\end{equation}
then $p_k > t_k$, $2|t_k$, and $C(p_k,q_k) \in \C_2$. Moreover,
\begin{equation}
\label{ap_1}
t_k = \sqrt{q_k^2-p_k^2} = \frac{q_k-1}{2},
\end{equation}
and so 
\begin{equation}
\label{ap_2}
\frac{1}{(2+\sqrt{3}) q_k} - \frac{2}{(2+\sqrt{3}) q_k^2} < \left| \frac{\sqrt{3}}{2} - \frac{p_k}{q_k} \right| < \frac{1}{2\sqrt{3}\ q_k} \longrightarrow 0, \text{ as } k \rightarrow \infty,
\end{equation}
and more precisely
\begin{equation}
\label{aphx1}
241 (7+4\sqrt{3})^{k-1} = 241 \times (13.928...)^{k-1} < q_k < 241 \times (13.947)^{k-1}.
\end{equation}
Hence, by (\ref{pack1}), for each such $C(p_k,q_k)$ the corresponding circle packing density is
\begin{equation}
\label{pack2}
\frac{\pi}{\sqrt{12}} \left( \frac{1}{1 + \frac{1}{723 (7+4\sqrt{3})^{k-1}}} \right) < \delta_{p_k,q_k} = \frac{\pi q_k}{4p_k} < \frac{\pi}{\sqrt{12}} \left( \frac{1}{1 + \frac{0.920...}{723 \times (13.947)^{k-1}}} \right),
\end{equation}
so $\delta_{p_k,q_k} \rightarrow \frac{\pi}{\sqrt{12}} = 0.9069... = \delta(\Lambda_h)$ as $k \rightarrow \infty$, and the quadratic form $Q_{p_k,q_k}(x,y)$ as in (\ref{Q_pq}) satisfies
\begin{eqnarray}
\label{ap_3}
\frac{1}{q_k} Q_{p_k,q_k}(x,y) & = & x^2 + \left( \frac{q_k-1}{q_k} \right) xy + y^2 \nonumber \\
& \longrightarrow & Q_h(x,y) := x^2 + xy + y^2, \text{ as } k \rightarrow \infty,
\end{eqnarray}
where $Q_h(x,y)$ is the norm form of $\Lambda_h$ with respect to the basis matrix as in (\ref{hex}).
\end{lem}

\proof
We start by proving (\ref{ap_1}). Let $(p_k,t_k,q_k)$ be given by (\ref{ap_0}), then
\begin{equation}
\label{app0}
\left( \begin{matrix} p_k \\ t_k \\ q_k \end{matrix} \right) = \left( \begin{matrix} 7&-4&8 \\ 4&-1&4 \\ 8&-4&9 \end{matrix} \right)^k \left( \begin{matrix} 15 \\ 8 \\ 17 \end{matrix} \right).
\end{equation}
We argue by induction on $k$. First notice that $p_1=209$, $t_1=120$, and $q_1=241$, so that $p_1 > t_1$, $2|t_1$, and (\ref{ap_1}) is satisfied. Now assume this holds for $(p_{k-1},t_{k-1},q_{k-1})$. By (\ref{app0}),
\begin{equation}
\label{app0.1}
\left( \begin{matrix} p_k \\ t_k \\ q_k \end{matrix} \right) = \left( \begin{matrix} 7p_{k-1}-4t_{k-1}+8q_{k-1} \\ 4p_{k-1}-t_{k-1}+4q_{k-1} \\ 8p_{k-1}-4t_{k-1}+9q_{k-1} \end{matrix} \right),
\end{equation}
and so
$$p_k = t_k + (3p_{k-1}-3t_{k-1}+4q_{k-1}) > t_k,$$
since $p_{k-1} > t_{k-1}$, as well as
$$t_k = 4(p_{k-1}+q_{k-1}) - t_{k-1}$$
is divisible by 2, since $2|t_{k-1}$, and finally
$$\frac{q_k-1}{2} = 4p_{k-1}-t_{k-1}+4q_{k-1} + \left( \frac{q_{k-1}-1}{2} - t_{k-1} \right) = t_k,$$
since $\frac{q_{k-1}-1}{2} = t_{k-1}$. The conclusion follows by induction.

Next we derive (\ref{ap_2}) from (\ref{ap_1}). Notice that by squaring both sides of (\ref{ap_1}) and rearranging terms, we immediately obtain
\begin{equation}
\label{app1}
\left( \frac{p_k}{q_k} - \frac{\sqrt{3}}{2} \right) \left( \frac{p_k}{q_k} + \frac{\sqrt{3}}{2} \right) = \frac{q_k-2}{2q_k^2},
\end{equation}
and since $\frac{p_k}{q_k} < 1$, we have
$$\left| \frac{p_k}{q_k} - \frac{\sqrt{3}}{2} \right| > \frac{q_k-2}{(2+\sqrt{3}) q_k^2},$$
which is the lower bound of (\ref{ap_2}). For the upper bound, we rewrite (\ref{app1}) as
$$\left| \frac{p_k}{q_k} - \frac{\sqrt{3}}{2} \right| = \frac{q_k-2}{q_k(2p_k+\sqrt{3}q_k)} < \frac{1}{2p_k+\sqrt{3}q_k} < \frac{1}{2\sqrt{3}\ q_k},$$
since $\sqrt{3}q_k < 2p_k$. It is also clear that $q_k \rightarrow \infty$ as $k \rightarrow \infty$.

To prove (\ref{aphx1}), we first notice that $q_1=241$. Moreover, by (\ref{ap_2}), the sequence $p_k/q_k$ is monotone decreasing and converges to $\sqrt{3}/2$, therefore
$$\frac{\sqrt{3}}{2} \leq \frac{p_k}{q_k} \leq \frac{p_1}{q_1} = \frac{209}{241},$$
for every $k \geq 1$. Then, by (\ref{app0.1}) and (\ref{ap_1}),
$$q_k=8p_{k-1}-4t_{k-1}+9q_{k-1} \geq (7+4\sqrt{3}) q_{k-1}+2 > (7+4\sqrt{3}) q_{k-1},$$
and 
$$q_k=8p_{k-1}-4t_{k-1}+9q_{k-1} \leq \left( 7 + \frac{8 \times 209}{241} \right) q_{k-1} + 2 < 13.947 \times q_{k-1}.$$
The inequalities (\ref{aphx1}) follow by induction on $k$.

To prove (\ref{pack2}), notice that upper bound (\ref{ap_2}) implies that
$$\frac{p_k}{q_k} < \frac{\sqrt{3}}{2} + \frac{1}{2\sqrt{3}q_k} = \frac{\sqrt{3}}{2} \left( 1 + \frac{1}{3q_k} \right) < \frac{\sqrt{3}}{2} \left( 1 + \frac{1}{723 (7+4\sqrt{3})^{k-1}} \right),$$
where the last inequality is obtained by applying by the lower bound of (\ref{aphx1}). Then the lower bound of (\ref{pack2}) follows. To obtain the upper bound of (\ref{pack2}), combine the lower bound of (\ref{ap_2}) with the upper bound of (\ref{aphx1}) in a similar manner.

Finally notice that (\ref{ap_3}) follows immediately from (\ref{ap_1}) and the fact that $q_k \rightarrow \infty$ as $k \rightarrow \infty$, and this completes the proof.
\endproof

\noindent
{\it Proof of Theorem \ref{app_hex}.}
Let $\left< \Lambda_{p_k,q_k} \right>$ be the sequence of similarity classes corresponding to the triples $(p_k,t_k,q_k)$ as defined in (\ref{ap_0}), then (\ref{ap_2}) guarantees convergence of this sequence to the similarity class $\left< \Lambda_h \right>$ with respect to the metric $d_s$ on $\Sim(\real^2)$, and also implies (\ref{aphx}), since $q_k \geq q_1 = 241$. The fact that $q_k=O(14^k)$ follows immediately from (\ref{aphx1}). To prove (\ref{aphx2}), assume that there exists some similarity class $\left< \Lambda_{p,q} \right> \neq \left< \Lambda_{1,1} \right>$ such that 
$$d_s(\Lambda_h,\Lambda_{p,q}) = \frac{p}{q} - \frac{\sqrt{3}}{2} \leq \frac{1}{3\sqrt{3}\ q},$$
which implies that $(3\sqrt{3}p-1)^2 \leq \frac{81}{4} q^2$, and therefore
$$q^2-p^2 \geq \left( \frac{1}{2} \right)^2 \frac{27q^2-24\sqrt{3}p+4}{27},$$
where
$$\frac{27q^2-24\sqrt{3}p+4}{27} > q^2-\frac{8p}{3\sqrt{3}} > q(q-\sqrt{3}) > 1,$$
since if $q > 1$, then $q \geq 13$. Hence
$$\sqrt{q^2-p^2} > \frac{1}{2},$$
which contradicts the fact that either $(p,\sqrt{q^2-p^2},q)$ or $(\sqrt{q^2-p^2},p,q)$ is in $\PP$, and so (\ref{aphx2}) must be true for each similarity class of the form $\left< \Lambda_{p,q} \right>$. This completes the proof of the theorem.
\boxed{ }
\smallskip

Finally, Corollary \ref{app_cor} follows immediately from (\ref{pack2}).
\smallskip

The approximation result of Theorem \ref{app_hex} is also interesting since the similarity class $\left< \Lambda_h \right>$ has a number of important properties: besides providing the optimal circle packing and minimizing Epstein zeta function, as mentioned in Section~1, it also solves the related minimization problem for the height of flat tori in dimension 2 (see \cite{chiu} for details), as well as the quantizer problem in dimension~2 (see \cite{conway} for details). Let us also recall that a lattice $\Lambda$ is called {\it perfect} if any real symmetric matrix $A$ in the corresponding dimension can be represented as
$$A = \sum_{\bx \in S(\Lambda)} \alpha_{\bx} \bx \bx^t,$$
where $S(\Lambda)$ is the set of minimal vectors of $\Lambda$ as in Section~1, each $\bx$ is written as a column vector, and each $\alpha_{\bx}$ is a real number. It is not difficult to see that for a lattice $\Lambda$ in $\real^2$ to be perfect, the cardinality of $S(\Lambda)$ must be six, meaning that the only perfect lattices in $\real^2$ come from $\left< \Lambda_h \right>$. Moreover, $\left< \Lambda_h \right>$ is {\it strongly perfect}, meaning that it supports a spherical 5-design: we say that a lattice $\Lambda$ in $\real^N$ (and hence its similarity class) supports a {\it spherical $t$-design} for $t \in \zed_{>0}$ if for every homogeneous polynomial $f(\bx)$ of degree $\leq t$ with real coefficients
\begin{equation}
\label{design}
\int_{\Ss^{N-1}} f(\bx) d\bx = \frac{1}{|S(\Lambda)|} \sum_{\bx \in S(\Lambda)} f(\bx),
\end{equation}
where $\Ss^{N-1}$ is the unit sphere in $\real^N$ with the canonical measure $d\bx$ on it, normalized so that $\int_{\Ss^{N-1}} d\bx = 1$. No other similarity class in $\Sim(\real^2)$ supports a spherical 5-design (or 4-design), and $\left< \Lambda_{1,1} \right>$ is the only other similarity class that supports a spherical 3-design (or 2-design); such similarity classes are called {\it strongly eutactic} (clearly, every lattice supports a 1-design). For detailed information on perfect and eutactic lattices see \cite{martinet}, especially chapter 16 for connections to spherical designs.
\bigskip

\section{Diophantine approximation by quotients of Pythagorean triples}

In this section we first prove Theorem \ref{wr_dense}. It follows immediately from the following direct consequence of a theorem of Hlawka \cite{hlawka} on simultaneous Diophantine approximation by quotients of Pythagorean triples, which we state here.

\begin{thm} \label{hlawka} Let $x \in (0,1)$ be a real number. Then there exist infinitely many Pythagorean triples $\left( p, \sqrt{q^2-p^2}, q \right)$ such that
\begin{equation}
\label{pyth_bound}
\left| x - \frac{p}{q} \right| \leq \frac{2 \sqrt{2}}{q}.
\end{equation}
\end{thm}

\noindent
{\it Proof of Theorem \ref{wr_dense}.}
Recall that
$$d_s( \Lambda, \Lambda_{p,q} )  = \left| \sin \theta(\Lambda) - \sin \theta(\Lambda_{p,q}) \right| = \left| \sin \theta(\Lambda) - \frac{p}{q} \right|,$$
and apply Theorem \ref{hlawka} with $x=\sin \theta(\Lambda)$.
\endproof

Moreover, we can say that the set $\{ \left< \Lambda_{p,q} \right> : p/q \in \s\}$ of similarity classes of WR sublattices of $\zed^2$ is equidistributed in the set $\Sim(\real^2)$ of similarity classes of all WR lattices in $\real^2$ in the following sense. It is a well known fact that the map 
$$t \mapsto \left( \frac{1-t^2}{1+t^2}, \frac{t}{t^2-1} \right)$$
is a bijection from the set of rational numbers onto the set of all rational points on the unit circle. Ordering $\que$ as the set of Farey fractions induces an ordering on the set of rational points on the unit circle, and hence on the set $\s$ of $y$-coordinates of such points that fall in the interval $\left[ \frac{\sqrt{3}}{2}, 1 \right]$. Now, it is a well known fact that Farey fractions are uniformly distributed ($\md 1$).
\smallskip

As a side remark, we can also use Theorem \ref{hlawka} to approximate points on a unit circle with rational points on the same circle. 

\begin{cor} \label{circle} Let $(x,y)$ be a point on the unit circle. Then either $x,y \in \{0, \pm 1 \}$, or there exist infinitely many rational points $(p/q, r/q)$ on the same circle such that
\begin{equation}
\label{circle_bound}
\max \left\{ \left| x - \frac{p}{q} \right|, \left| y - \frac{r}{q} \right| \right\} \leq \frac{2 \sqrt{2}}{q}.
\end{equation}
\end{cor}

\proof
First notice that it suffices to prove the statement of this corollary for the case $0 < x,y < 1$, namely the case when the point in question lies in the first quadrant, since any other point on the circle can be obtained from those in the first quadrant by a rational rotation. Let $c$ be an arbitrary real number in the interval $(0,1)$, then either 
\begin{equation}
\label{xy1}
0 < x \leq \sqrt{1-c^2} < 1,\ c \leq y < 1,
\end{equation}
or
\begin{equation}
\label{xy2}
0 < y \leq \sqrt{1-c^2} < 1,\ c \leq x < 1.
\end{equation}
First assume that (\ref{xy1}) holds. By Theorem \ref{hlawka}, there exist infinitely many Pythagorean triples $\left( p, r, q \right)$ with $r = \sqrt{q^2-p^2}$ which satisfy (\ref{pyth_bound}). Then:
\begin{eqnarray}
\label{c1}
\frac{2 \sqrt{2}}{q} \geq \left| x - \frac{p}{q} \right| & = & \left| \sqrt{1-y^2} - \sqrt{ 1-\frac{r^2}{q^2}} \right| = \frac{ \left| \frac{r^2}{q^2} - y^2 \right| }{\sqrt{1-y^2} + \sqrt{1 - \frac{r^2}{q^2}}}  \nonumber \\
& = & \frac{ \frac{r}{q} + y }{\sqrt{1-y^2} + \sqrt{1 - \frac{r^2}{q^2}}} \left| y - \frac{r}{q} \right| \geq \frac{c \left( 1 + \frac{n}{n+1} \right)}{2 \sqrt{1- \frac{n^2}{(n+1)^2} c^2}} \left| y - \frac{r}{q} \right|.
\end{eqnarray}
The last inequality is true because $\frac{w+z}{\sqrt{1-w^2}+\sqrt{1-z^2}}$ is an increasing function in both variables for $0 < z,w < 1$; since $y \geq c$, we can pick $q$ large enough so that $r/q$ would have to be sufficiently close to $y$ so that $r/q \geq \frac{n}{n+1} c$ for some $n \in \zed_{>0}$, then $r/q + y \geq c \left( 1 + \frac{n}{n+1} \right)$, and $\sqrt{1-y^2} + \sqrt{1 - \frac{r^2}{q^2}} \leq 2 \sqrt{1- \frac{n^2}{(n+1)^2} c^2}$. Then (\ref{c1}) implies:
\begin{equation}
\label{c2}
\left| y - \frac{r}{q} \right| \leq \frac{\sqrt{1- \frac{n^2}{(n+1)^2} c^2}}{c \left( 1 + \frac{n}{n+1} \right)} \times \frac{4 \sqrt{2}}{q}.
\end{equation}
Since our choice of $c \in (0,1)$ and positive integer $n$ was arbitrary, we can for instance choose 
\begin{equation}
\label{c3}
c = \frac{2n+2}{\sqrt{8n^2+4n+1}},
\end{equation}
and take $n=2$, in which case, combining (\ref{pyth_bound}), (\ref{c2}), and (\ref{c3}), we obtain (\ref{circle_bound}).

If, on the other hand, (\ref{xy2}) holds instead of (\ref{xy1}), simply repeat the above argument interchanging $x$ with $y$ and $p/q$ with $r/q$. This completes the proof.
\endproof

A related result has also been obtained by Kopetzky in  \cite{kopetzky1} (also see \cite{kopetzky2}), however his bounds are different in flavor in the sense that the constants in the upper bounds depend on $x$ and $y$. Notice that the bound of Corollary \ref{circle} can be easily extended to any rational ellipse.

\begin{cor} \label{ellipse} Let $(x,y)$ be a point on the ellipse $E$, given by the equation
$$\left( \frac{x}{a} \right)^2 + \left( \frac{y}{b} \right)^2 = 1,$$
where $a,b$ are positive rational numbers. Then either $(x,y) = (\pm a, 0), (0, \pm b)$, or there exist infinitely many rational points $(p/q, r/q)$ on the same ellipse such that
\begin{equation}
\label{ellipse_bound}
\max \left\{ \left| x - \frac{p}{q} \right|, \left| y - \frac{r}{q} \right| \right\} \leq \frac{2 \sqrt{2} \max\{ a,b \}}{q}.
\end{equation}
\end{cor}

\proof
Notice that the map $(x,y) \mapsto (x/a,y/b)$ is a bijection between $E$ and the unit circle, which takes rational points to rational points. Now apply Corollary \ref{circle} to points of the form $(x/a,y/b)$.
\endproof
\bigskip

{\bf Acknowledgment.} I would like to thank Pavel Guerzhoy and the referees for their helpful comments on the subject of this paper. I would also like to acknowledge the wonderful hospitality of Institut des Hautes \'{E}tudes Scientifiques in Bures-sur-Yvette, France, where a part of this work has been done.

\bibliographystyle{plain}  
\bibliography{wr_sim}   

\end{document}